\documentclass[12pt]{article}

\usepackage{amsmath,amsthm,amssymb}

\begin{document}

\theoremstyle{plain}

\newtheorem{Thm}{Theorem}[section]
\newtheorem{Prop}{Proposition}[section]
\newtheorem{Defi}{Definition}[section]
\newtheorem{Cor}{Corollary}[section]
\newtheorem{Lem}{Lemma}[section]

\makeatletter\renewcommand{\theequation}{%
\thesection.\arabic{equation}}
\@addtoreset{equation}{section}\makeatother

\setlength{\baselineskip}{16pt}
\newcommand{\vlimsup}{\mathop{\overline{\lim}}}
\newcommand{\vliminf}{\mathop{\underline{\lim}}}
\newcommand{\Av}{\mathop{\mbox{Av}}}
\newcommand{\spec}{{\rm spec}}
\def\textmc{\rm}
\def\({(\!(}
\def\){)\!)}
\def\R{\mathbb{R}}
\def\Z{\mathbb{Z}}
\def\N{\mathbb{N}}
\def\C{\mathbb{C}}
\def\T{\mathbb{T}}
\def\E{{\bf E}}
\def\H{\mathbb{H}}
\def\Prob{{\bf P}}
\def\M{{\cal M}}     
\def\F{{\cal F}}
\def\G{{\cal G}}
\def\D{{\cal D}}
\def\X{{\cal X}}
\def\A{{\cal A}}
\def\B{{\cal B}}
\def\L{{\cal L}}
\def\a{\alpha}
\def\b{\beta}
\def\e{\varepsilon}
\def\de{\delta}
\def\ga{\gamma}
\def\k{\kappa}
\def\la{\lambda}
\def\fa{\varphi}
\def\th{\theta}
\def\si{\sigma}
\def\t{\tau}
\def\om{\omega}
\def\De{\Delta}
\def\Ga{\Gamma}
\def\La{\Lambda}
\def\Om{\Omega}
\def\Th{\Theta}
\def\lan{\langle}
\def\ran{\rangle}
\def\lbr{\left(}
\def\rbr{\right)}
\def\const{\;\operatorname{const}} 
\def\dist{\operatorname{dist}} 
\def\Tr{\operatorname{Tr}}
\def\quadd{\qquad\qquad}
\def\n{\noindent}
\def\beq{\begin{eqnarray*}}
\def\eeq{\end{eqnarray*}}
\def\supp{\mbox{supp}}
\def\beqn{\begin{equation}}
\def\eeqn{\end{equation}}
\def\bp{{\bf p}}
\def\sg{{\rm sgn\,}}
\def\1{{\bf 1}}
\def\pf{{\it Proof.}}
\def\v2{\vskip2mm}
\def\n{\noindent}
\def\z{{\bf z}}
\def\x{{\bf x}}
\def\y{{\bf y}}
\def\bv{{\bf v}}
\def\be{{\bf e}}
\def\pr{{\rm pr}}

\def\cp{{\rm Cap}}
\def\tst12{{\textstyle\frac12}}

\begin{center}
{\Large Density of space-time distribution\\
 of Brownian first hitting
 of a disc and a ball.  } \\
\vskip4mm
{K\^ohei UCHIYAMA} \\
\vskip2mm
{Department of Mathematics, Tokyo Institute of Technology} \\
{Oh-okayama, Meguro Tokyo 152-8551\\
e-mail: \,uchiyama@math.titech.ac.jp}
\end{center}

\vskip18mm
\n
{\it running head}:   space-time distribution of Brownian first hitting

\vskip2mm
\n
{\it key words}: harmonic measure for heat operator, Brownian hitting time, caloric measure, parabolic measure, space-time distribution
\vskip2mm
\n
{\it AMS Subject classification (2010)}: Primary 60J65,  Secondary  60J45, 60J60. 

\vskip16mm

\begin{abstract}
We compute  the joint distribution of the site  and the time  at which a $d$-dimensional standard Brownian motion $B_t$ hits the surface of the ball  $ U(a) =\{|{\bf x}|<a\}$  for the first time.
The asymptotic form of its density is
obtained when either the hitting time or the starting site $B_0$ becomes large.  Our results  entail that if Brownian motion is started at ${\bf x}$ and conditioned to hit $U(a)$ at time $t$ for the first time,  the distribution of  the hitting site  approaches  the uniform distribution  or  the point mass at $a{\bf x}/|{\bf x}|$ according as  $|{\bf x}|/t$ tends to zero or infinity;  in each case  we provide a  precise asymptotic  estimate of the density.  In the case when  $|{\bf x}|/t$ tends to a positive constant we show the convergence of the density  and  derive an analytic  expression of  the limit density.
   \end{abstract}
\vskip6mm

\n
\section{ Introduction}

The harmonic measure (also called caloric measure or caloric measure in the present context \cite{W1}) of the  unbounded space-time domain
$$D=\{(\x,t)\in \R^d\times (0,\infty) :  |\x|>a\}$$ $(a>0)$
for the heat operator $\frac12\De -\partial_t$ consists of two components, one supported by  the initial  time boundary $t=0$  and the other by the lateral boundary $\{|\x|= a\}\times \{t>0\}$.  The former  one is nothing but the measure whose density is given by the heat kernel for physical space $|\x|>a$ with Dirichlet zero boundary condition. This paper concerns  the latter,   aiming  to find   a precise  asymptotic form of it  when  the distance of the reference point from the boundary 
becomes large.   In the probabilistic term this latter part  is   given by  the joint distribution, $H(\x, dtd\xi)$,  of the site  $\xi$ and the time  $t$ at which the $d$-dimensional standard Brownian motion hits the surface of the ball  $ U(a) =\{|\x|<a\}$  for the first time: given a bounded continuous function $\fa(\xi,t)$ on the lateral boundary of  $D$, the bounded solution $u=u(\x,t)$ of the heat equation 
$(\frac12 \De -\partial_t)u=0$ in $D$ satisfying the  boundary condition
$$u(\xi, t) = \fa(\xi, t)  \,\,\, (|\xi|=a, t>0) \quad\mbox{and} \quad u(\x,0)= 0 \,\, \, (|\x|>a)$$
can be expressed in the boundary integral 
$$u(\x,t) =\int_0^t \int_{|\xi|=a}\fa(\xi, t-s) H(\x, dsd\xi).$$
The probability  measure $H(\x, dtd\xi)$  has a smooth  density, which may be factored into the product of    the hitting time density  and the density for the hitting site distribution conditional on the hitting  time.
While the asymptotic forms of the first factor are  computed in several recent papers \cite{Ubh},  \cite{BMR}, \cite{HM}, \cite{Ubes}, the latter seems to be rarely investigated and   in this paper we carry out   the computation of it.  Consider  the hitting  site   distribution  of $\partial U(a)$ for the Brownian motion  conditioned to start at  $\x \notin U(a)$ and hit  $U(a)$ at time $t$ for the first time.
 It would be intuitively clear that the conditional distribution of the hitting site  becomes nearly uniform   on the sphere for large  $t$  if $|\x|$ is small relative to $t$, while one may speculate that it  concentrates about the point $a\x/|\x|$ as $|\x|$ becomes very large in comparison with $t$. Our results  entail that in the limit there appears  the uniform distribution   or  the point mass at $a\x/|\x|\in \partial U(a)$ according as  $|\x|/t$ tends to zero or infinity;  in each case  we provide  a certain  exact estimate of the density.  In the case when  $|\x|/t$ tends to  a positive constant the conditional distribution has a limit, of which  we   derive an analytic  expression for  the   density.  Using  these results together with the estimates of hitting time density  obtained in \cite{Ubes}  we can compute the hitting  distributions of bounded Borel sets, as is carried out in a separate paper \cite{Ucarl_B}.  When $|\x|/t$ tends to become large,  the problem is comparable to that for the hitting distribution for the Brownian motion   with a large constant drift started at $\x$ and for the latter process  one may  expect that the distribution is uniform if it is projected on  the cross section of  $U(a)$ cut with  the plane  passing through the origin and perpendicular to the unit vector  $\x/|\x|$. This is   true in the sense of weak  convergence of  measures,  but  in a finer measure the distribution is not  flat:  the density of the projected distribution  has large values along the circumference of the cross section.  For such computation it is crucial  to have a certain delicate estimate of the hitting  distribution  for $t$ small, which  we  also  provide in this paper.

\section{Notation and Main Results}
In this section  we present main results obtained in this paper, of which some detailed statements may be given later sections.  Before doing that, we give basic notation used throughout   and state the results on the hitting time distribution from \cite{Ubes}.
\v2\n
\vskip20mm
{\sc  {\bf 2.1.} Notation.} \,

We fix the radius $a>0$  of the Euclidian  ball  $U(a)= \{\x\in \R^d: |\x|<a\}$  ($d=2,3,\ldots$). Let $P_\x$ be the probability law of a $d$-dimensional  standard Brownian motion, denoted by $B_t, t\geq 0$,   started at $\x\in \R^d$ and $E_\x$ the expectation under $P_\x$. We usually write $P$ and $E$ for $P_{\bf 0}$ and   $E_{{\bf 0}}$, respectively, where ${\bf 0}$ designates the origin of $\R^d$.

The following notation is used throughout the paper.  
 \beq
&&\nu =\frac{d}{2}-1 \quad (d=1,2,\ldots);\\ 
&&{\bf e} = (1, 0,\ldots,0) \in \R^d; \\
&&\sigma_a =\inf \{t>0: |B_t|\leq a\}; \\
&&q_a^{(d)}(x,t) =\frac{d}{dt}P_\x[\sigma_a \leq t] \quad (x=|\x|>a).
\\
&& p_t^{(d)}(x)   = (2\pi t)^{-d/2} e^{-x^2/2t}.\\
&& \La_\nu(y)= \frac{(2\pi)^{\, \nu+1}}{2y^{\nu}K_{\nu}(y)}\quad (y>0).\\
 &&\om_{d-1} =2\pi^{d/2}/ \Ga(d/2) \, \,\,(\mbox{the area of $d-1$ dimensional unit sphere}).\\
 &&\mu_d = \om_{d-1}/\om_{d-2} = \sqrt{\pi}\,  \Ga(\nu+\tst12)/\Ga(\nu+1).
 \eeq
Here $K_\nu$ is the modified Bessel function of the  second kind of order $\nu$.  We usually write 
$x$ for $|\x|$, $\x\in \R^d$  (as above);  $d=2\nu +2$ and $\nu$ are used interchangeably;  and we sometime write  $q^\nu(x,t)$ for  $q^{(d)}(x,t)$   when  doing so gives rise to no confusion and  facilitates computation or exposition and also   $B(t)$ for  $B_t$ for typographical reason.  When working on the plane we often tacitly  use complex notation to denote points of it, for instance a point of the circle $\partial U(a)$ is indicated as $ae^{i\th}$ with  $\th$  denoting the (well-defined) argument of the point.

We write $x\vee y$ and $x\wedge y$ for the maximum and  minimum of real numbers $x, y$, respectively; $f(t) \sim g(t)$ if
 $f(t)/g(t)\to 1$ in any  process of taking limit.   The symbols  $C, C_1, C',$  etc,  denote universal positive  constants whose precise values are unimportant; the same symbol may takes different values  in different occurrences. 
 
\v2\n
{\sc {\bf 2.2.}  Density of Hitting Time  Distribution.} 

Here we state  the  results from \cite{Ubes}  on  $q_a^{(d)}(x,t)$, the density    for  $\sigma_a$.
The definition of   $q^{(d)}(x,t)$    may be naturally extended to  Bessel processes of order $\nu$ and the results concerning it given below may be applied to such extension if $\nu\geq 0$. 
\v2\n
{\bf Theorem A.} 
{\it Uniformly for $x > a$, as $t\to\infty$, 
\beqn\label{R2}
q_a^{(d)}(x,t)  \,\sim\, a^{2\nu}\La_\nu\bigg(\frac{a x}{t}\bigg) p^{(d)}_t(x) \bigg[1-\bigg(\frac{a}{x}\bigg)^{2\nu}\bigg] \qquad  (d\geq 3)
\eeqn
and  for $d=2$,}
\beqn\label{R21}
q_a^{(2)}(x,t)  = p^{(2)}_t(x) 
\times
 \left\{ \begin{array} {ll}  
  {\displaystyle \frac{4\pi\lg(x/a)\,}{(\lg (t/a^2))^2}\Big(1+o(1)\Big)    } \quad& (x  \leq \sqrt t\,),\\ [5mm]
 {\displaystyle  \La_0\bigg(\frac{a x}{t}\bigg) \Big(1+o(1)\Big)  }\quad&(  x   > \sqrt t\,).
  \end{array} \right.
\eeqn

\v2
 From the known properties of $K_\nu(z)$ it follows that
 \beqn\label{lambda}
\La_\nu(y) = (2\pi)^{\nu+1/2} y^{-\nu+1/2}\, e^{y} ( 1+O(1/y)) \quad \mbox{as}\quad y\to\infty;
\eeqn
\[
\La_\nu(0) = \frac{2\pi^{\nu+1}} {\Ga(\nu)} ( = \nu\om_{d-1})
\quad\mbox{for}\quad\nu >0;\quad \La_0(y) \sim \frac{\pi}{-\lg y}  \quad \mbox{as}\quad y \downarrow 0.
\]

\v2\n
{\bf Theorem B.}
\, 
{\it For each $\nu\geq 0$ it holds that  uniformly
 for  all $t>0$ and $x>a$, }
\beqn\label{result6}
 q_a^{(d)}(x,t) =  \frac{x-a}{\sqrt{2\pi}\, t^{3/2}} e^{-(x-a)^2/2t}\bigg(\frac{a}{x}\bigg)^{(d-1)/2}\Bigg[1+ O\bigg(  \frac{t}{ax}\bigg)\Bigg].
 \eeqn

 \v2\n
 {\sc Remark 1.}\,   Under certain constraints   on $x$ and $ t$  some finer  error estimates in the formulae of Theorem A  are given in \cite{Ubh} ($d=2$, $|x|<\sqrt t$))   and in \cite{Ubes} ($|x|/t\to\infty$).   The formula (\ref{result6})  of Theorem B   is sharp only if  $x/t \to\infty$. The case $t\to\infty$ of it  is contained in Theorem A apart from the error estimate. A better error estimate is obtained   in  \cite{BMR} by  a purely  analytic approach.  A probabilistic proof of (\ref{result6})  is found  in \cite{Usaus2}. We  shall use  (\ref{result6}) primarily  for the case  $0<t<a^2$.

 \v2\n
{\sc Remark 2} (Scaling property). From the scaling property of Brownian motion  it follows that 
\[
 q^{(d)}_a(x,t) = a^{-2} q^{(d)}_1 (x/a, t/a^2); \mbox{and}
 \]
\[ \frac{P_{x{\bf e}}[B(\sigma_a)\in d\xi\,|\, \sigma_a =t] }{m_a(d\xi)} = \frac{P_{(x/a){\bf e}}[B(\sigma_1)\in d\xi'\,|\, \sigma_1 =t/a^2] }{m_1(d\xi')}\Big|_{\xi' = \xi/a}
 \]
for all dimensions $\geq 2$. 
Even though because of this we can obtain the result for $a\neq 1$ by simply substituting  $t/a^2$
and $x/a$ in place of  $t$ and $x$, respectively,  in the formula for $a=1$, in the above we have exhibited the formula  for $q^{(d)}_a(x,t)$ with $a>0$ arbitrary.  We shall follow  this example  in  stating the results of the present work. It is warned that  we are not so
scrupulous in doing that: in particular,   to indicate the constrains of  $t$ (and/or $x$)  we often  simply write $t>1$ when we should  write $t>a^2$ for instance.
 

\v2\n
{\sc {\bf 2.3.} Density  of Hitting Site Distribution Conditional on $\sigma_a =t$.}  \, 

For finding 
  the asymptotic form of  the hitting distribution, with that  of $q^{(d)}(x,t) $ being given in {\bf 2.2},   it remains to estimate  the conditional density $P_{\x}[B_t\in d\xi | \sigma_a=t]/d\xi$.    Before stating the results on it we shall consider  the argument of the hitting site $B(\sigma_a)$ in the case $d=2$, when the winding number around the origin is naturally  associated with  the process.

 \v2
{\sc {\bf 2.3.1.}  Density for $\arg B(\sigma_a)$ (Case $d=2$)}. \, Let $\arg B_t\in \R$ be the argument of $B_t$   (regarded as  a complex Brownian motion),  which  is a.s. uniquely determined by continuity under the  convention  $\arg B_0\in (-\pi,\pi]$.  The function of $\la\in \R$  defined by
\beqn\label{Def}
 \Phi_a(\la;v)  =\frac{K_0(av)}{K_{\la}(av)} \quad  (v>0).
 \eeqn
turns out to be a characteristic function of a probability distribution on $\R$.
\begin{Thm}\label{prop2} \quad
${\displaystyle \Phi_a(\la;v) = \lim_{x/t\to v} E_{x{\bf e}}[e^{i\la\arg B(\sigma_a)}| \sigma_a = t].}  $
\end{Thm}
Since  $\Phi_a(\la;v)$  is continuous at $\la=0$, Theorem \ref{prop2} shows  that  the conditional law of $\arg B(\sigma_a)$ converges to the probability law whose   characteristic function  equals  $\Phi_a(\la;v)$.  
In fact     $\Phi_a$   is smooth in $\la$, so that  the limit  law has a density. If $f_a(\th;v)$ denotes the density,  then
$$\Phi_a(\la;v) = \int_{-\infty}^\infty  e^{i\la \th }f_a(\th;v)d\th \qquad (\la \in \R);$$
 we shall see that the density of the conditional law converges:
\begin{eqnarray}\label{f} 
f_a(\th;v) = \lim_{x/t\to v}  \frac{P_{x{\bf e}}[\arg B_t\in d\th| \sigma_a = t] }{d\th} \qquad\quad  (-\infty <\th<\infty, v>0)
\end{eqnarray}
 (Section {\bf 3.2.2}).
By (\ref{Def}) 
$$\Phi_{a}(\la;0+) =0 \quad  (\la\neq 0) \quad \mbox{ and} \quad 
 \Phi_{a}(\la;+\infty) =1,$$
which shows  that the probability   $f_a(\th;v)d\th$ concentrates  in the limit at infinity  as $v\downarrow 0$ and
 at zero   as $v\to\infty$.
 Since  for $0<y<\infty$,
\beqn\label{K/nu}
\lg K_\la(y)\sim |\la|\lg |\la|  \quad\mbox{as}\quad \la \to\pm \infty,
\eeqn
for each $v>0$,  $f_a(\,\cdot\,;v)$ can be extended to an entire function;  in particular its support (as a function on $\R$)  is the whole real line and we can then readily infer that $f_a(\th;v)>0$ for all $\th$ (see  (\ref{f_v})).   $K_{i\eta}(av)$ is an entire function of $\eta$ and has zeros on and only on the real axis.  If $\eta_0$ is its smallest positive zero, then 
\v2
\quad\qquad  $\int_0^\infty f_a(\th;v)e^{\eta \th}d\th $ is finite  or infinite  according as $\eta<\eta_0$ or $\eta\geq \eta_0$;
\v2\n
 it  can be  shown  that $0< \eta_0 - av\leq C(av)^{1/3}$ for $v>1$ and $-C(\lg av)^2 \leq \eta_0-\pi/|\lg av|<0$ for $v<1/2$ \cite{Uw}. 
 
 The next result is derived independently of Theorem \ref{prop2} and in a quite different way.
\begin{Prop}\label{thm0} \,  For $v>0$
\beqn\label{EQ0}
f_a(\th;v) \geq  \pi^{-1} a vK_0(av) \,e^{av\cos \th} \cos \th \qquad (|\th|\leq \tst12\pi).
\eeqn
\end{Prop}

\v2
{\sc {\bf 2.3.2.}\,  Density for Hitting Site}. \, 
 Let $m_a(d\xi)$ denote the uniform  probability distribution on $\partial U(a)$, namely 
 $m_a(d\xi) = (\omega_{d-1}a^{d-1})^{-1}|d\xi|$, where  $\omega_{d-1}$ denotes  the area of the  $d-1$ dimensional unit sphere $\partial U(1)$,  $d\xi \subset \partial U(a)$ an surface element and  $|d\xi|$ its Lebesgue measure. Let  ${\rm Arg}\,  z$, $z\in \R^2$ denote the principal value $\in (-\pi,\pi]$  of $\arg z$.

\begin{Thm}\label{thm1.1}  \,\,  {\rm (i)} If $d=2$, uniformly for $\th\in (-\pi, \pi]$, as $ x/t\to 0$ and $t\to \infty$
$$\frac{P_{x{\bf e}}[{\rm Arg} \,B_t\in d\th  \,|\, \sigma_a =t] }{d\th} = \frac{1}{2\pi} + O\Big(\frac{x}{t} \ell(x,t)\Big),$$
where $\ell(x,t)= (\lg t)^{2}/ \lg (x+2a)$ if $a <x<\sqrt t$ and  $= \lg (t/x)$  if $x>\sqrt t$.

 {\rm (ii)} If $d\geq3$, uniformly
for $\xi\in \partial U(a)$, as $v= x/t\to 0$ and $t\to \infty$, 
$$\frac{P_{x{\bf e}}[B_t\in d\xi\,|\, \sigma_a =t] }{m_a(d\xi)} = 1+O\Big(\frac{x}{t}\Big).$$
\end{Thm}
\v2
The orders of magnitude for the error terms given  in Theorem \ref{thm1.1} are correct  ones  (see Theorem \ref{thm-1} and Corollary \ref{cor-12}).

Let $\th=\th(\xi)\in [0,\pi]$ denote the colatitude of a point $\xi\in \partial U(a)$ with $a\be$ taken to be the north pole, namely  $a\cos \th = \xi\cdot \be$.
\v2

\begin{Thm}\label{thm1.21}    For each $M>1$, uniformly for 
$0< v< M$ and $\xi \in \partial U(a)$, as $t\to\infty$ and $x/t \to v$
\beqn\label{lim}
\frac{P_{x{\bf e}}[B_t\in d\xi\,|\, \sigma_a =t] }{m_a(d\xi)} \,  \longrightarrow \, \sum_{n=0}^\infty \frac{K_\nu(av)}{K_{\nu+n}(av)} H_n(\th).
\eeqn
Here $\th =\th(\xi)$; and $H_0(\th) \equiv 1$ and for $n\geq 1$,
$$H_n(\th) = \left\{\begin{array} {ll} 2 \cos n \th   \quad &\mbox{if} \quad d=2,\\
(1+\nu^{-1}n) C_n^\nu(\cos \th) \quad 
  &\mbox{if} \quad d \geq 3,
 \end{array}\right.$$ 
where 
  $C_n^\nu(z)$ is  the Gegenbauer polynomial of order $n$ associated with $\nu$.
\end{Thm}

According to (\ref{K/nu}) the convergence of the series  appearing as the limit  in  (\ref{lim}) is quite fast.    For $d=2$, as one may notice,  (\ref{lim})  is obtained from Theorem \ref{prop2} by using
 Poisson summation formula.  The limit function represented by  the series approaches unity as $v\downarrow 0$ (uniformly in $\th$), so that the asserted uniformity of convergence implies that the density on the left converges to unity as $x/t\to 0$, comforming to Theorem \ref{thm1.1}.
\begin{Thm}\label{thm1.2}   \,  Uniformly for  $t>1$,  as $v:= x/t\to\infty$ 
\beq
&&\frac{P_{x{\bf e}}[B_t\in d\xi\,|\, \sigma_a =t] }{\om_{d-1}m_a(d\xi)}  \\
 &&\quad = \bigg(\frac{av}{2\pi}\bigg)^{(d-1)/2}e^{-av(1-\cos \th)}\cos \th \bigg[1+ O\Big(\frac{1}{av\cos^3 \th}\Big)\bigg]  
    \quad  \mbox{if} \quad  0 \leq \th \leq \frac12\pi - \frac1{(av)^{1/3}}, \\
&&\quad \asymp   \bigg(\frac{av}{2\pi}\bigg)^{(d-1)/2} e^{-av(1-\cos \th)}  \frac{1 }{ (av)^{1/3} }   \qquad  \mbox{if} \quad  \frac12\pi - \frac1{(av)^{1/3}} <  \th \leq \frac{\pi}{2} + \frac1{(av)^{1/3}},
\eeq
where $\th =\th(\xi)$; $f(t)\asymp g(t)$ signifies that $f(t)/g(t)$ is bounded away from zero and infinity.
\end{Thm}

Combined with   Theorem B,    Theorem \ref{thm1.2}  yields  an asymptotic result of the joint distribution of $(B_{\sigma_a}, \sigma_a)$. On noting  that
$(\frac{y}{2\pi})^{(d-1)/2}=[ye^{y}/\La_{\nu}(y)](1+ O(1/y))$ ($y>1)$,
 $\cos \th = \x\cdot \xi/ax$, 
$$e^{- av (1-\cos \th)} p_t^{(d)}(x-a) =  p^{(d)}_t(|\x-\xi|)$$
and $\cos \th \sim \frac12 \pi -\th$ as $\th \to\frac12 \pi$, we state  the first half of it  as the following 
\begin{Cor}\label{thm3.02} \,  Uniformly under the constraint   $\x\cdot \xi /ax > (av)^{-1/3}$ and $t>a^2$,  as $v:= x/t\to\infty$ 
$$
\frac{P_{x{\bf e}}[B(\sigma_a)\in d\xi,\, \sigma_a \in dt] }{\om_{d-1}m_a(d\xi)dt}  
= a^{2\nu}\frac{\x\cdot \xi}{t} p^{(d)}_t(|\x-\xi|)\Bigg[1+O\bigg( \frac1{av\cos ^3\th}\bigg)\Bigg].
$$
\end{Cor}

\v2 
 As is clear from Theorem \ref{thm1.2} the distribution of $B(\sigma_a)$ converges to the Dirac delta measure at $a\be$,  the north pole of $\partial U(a)$, as $v\to\infty$.    The distribution may be normalized so as to approach a positive multiple of the  non-degenerate  measure $\cos \th \, m_a(d\xi)$  in obvious manner,  even though the density has singularity along the circumference.
 The next corollary  states this  in terms  of the colatitude $\Th(\sigma_a) :=\th(B(\sigma_a)$ of $B(\sigma_a)$ (see also Lemma \ref{wc}).
 
\begin{Cor}\label{thm3.01}  As $v:= x/t\to\infty$ under  $t>a^2$
\beq
&&\bigg(\frac{2\pi}{av}\bigg)^{(d-1)/2}e^{av(1-\cos \th)} P_{x{\bf e}}[\Th(\sigma_a)\in d\th\,|\, \sigma_a =t] \\
&&\quadd\quad \Longrightarrow \quad \,\om_{d-2}\1(0\leq\th\leq \tst12 \pi) \cos \th\,\sin^{d-2}\th\, {d\th},
\eeq
where $\1({\cal S})$ is the indicator function of a statement ${\cal S}$, 
\lq \,$\Rightarrow$\rq \, signifies  the weak convergence of finite measures on $\R$  (in fact the convergence holds in the total variation norm) and $\om_0=2$.
\end{Cor}

The essential content involved in  Theorem \ref{thm1.2}  concerns the   two-dimensional Brownian motion even if it includes the higher-dimensional one (cf. Section 6).

The rest of the paper is organized as follows. In Section 3 we deal with the case when $x/t$ is bounded and prove Theorems \ref{prop2} through   \ref{thm1.21}. In Section 4 we provide several preliminary estimates of the hitting distribution density mainly for $t<1$, that prepare for verification  of Theorem \ref{thm1.2}  made in Section 5 for the case $d=2$ and in Section 6 for the case $d\geq 3$.  Proposition \ref{thm0} is obtained in Section 5.1 as a byproduct of a  preliminary result for the proof of Theorem \ref{thm1.2}. In Section 7 the results obtained  are applied   to  the corresponding problem for Brownian motion with drift. In the final section,  Appendix, we present a classical  formula for the hitting  distribution of $U(a)$ and give a comment on  an approach to the present  problem based on it.

\section{ Proofs of Theorems \ref{prop2} through   \ref{thm1.21} }

This section consists of three subsections.   In the first subsection  we let  $d=2$ and  prove  Theorem \ref{prop2}. The proofs of  Theorems \ref{thm1.1}  and \ref{thm1.21} are given  in the rest. The essential ideas for all of them  are already found in the first subsection.

Our proofs involve  Bessel processes of varying order $\nu$ and it is convenient to introduce notation specific to them.   Let   $X_t$ be  a Bessel process  of order $\nu\in \R$  and  denote   by $P_x^{BS(\nu)}$ and $E_x^{BS(\nu)}$ the probability law   of $(X_t)_{t\geq 0}$  started at  $x\geq 0$ and  the expectation w.r.t.  it, respectively.   If $\nu= -1/2$, it is a standard Brownian motion and we write $P^{BM}_x$ for $P_x^{BS(-1/2)}$. With this convention we 
 suppose $\nu\geq 0$ in what follows, so that $X_t\geq 0$ a.s.  under $P_x^{BS(\nu)}$ ($x\geq 0$). 
The  expression   $2\nu +2$ which is not integral may appear, while the letter $d$ always designates  a positive integer  signifying  the dimension of   $B_t$,  a $d$-dimensional Brownian motion under a probability law $P_\x$.  

  Let  $T_a$  denote the first passage time of $a$ for $X_t$:  $T_a = \inf\{t\geq 0: X_t =a\}$.

\v2\n
   {\sc {\bf 3.1.}   The characteristic function  of  $\arg B(\sigma_a)$ ($d=2$).}\, 
   
The proofs of Theorems \ref{prop2} and the case  $d=2$ of Theorems  \ref {thm1.1} and  \ref{thm1.21}   rest on  the following
\begin{Prop}\label{lem3.20}\,  For  $\la\in \R$, $x>a$ and $t>0$, 
 $$ E_{x{\bf e}}[ e^{i\la\arg B(\sigma_a)}\,|\,\sigma_a = t] 
 = \frac{q_a^{(2|\la|+2)}(x,t)}{q_a^{(2)}(x,t)} \bigg(\frac{x}{a}\bigg)^{|\la|}. 
 $$
\end{Prop}
 \n
 
 In this subsection we first exhibit   how this proposition leads  to  Theorem \ref{prop2},  then   
 prove  two lemmas  concerning  Bessel processes  and used  in later subsections as well,   and finally prove Proposition \ref{lem3.20} by using  these lemmas.

 \v2
 {\sc {\bf 3.1.1.}  Deduction  of Theorem \ref{prop2} from Proposition \ref{lem3.20}. }   On using  Theorem A and (\ref{lambda}) in turn,  as $x/t\to v>0$
\begin{eqnarray}
\frac{q_a^{(2|\la|+2)}(x,t)}{q_a^{(2)}(x,t)} \bigg(\frac{x}{a}\bigg)^{|\la|}
 &\sim & \bigg(\frac{x}{a}\bigg)^{|\la|}\frac{ a^{2|\la|}\La_{|\la|}(av)p^{(2|\la|+2)}_t(x) }{\La_0(av)p^{(2)}_t(x)}  
 \nonumber\\
&\sim&\frac{K_{0}(av)}{K_{|\la|}(av)}. \label{K/K}
\end{eqnarray}
Noting that $K_{-\nu} (z)=K_\nu(z)$,   we obtain the identity of Theorem  \ref{prop2} according to Proposition \ref{lem3.20}. \qed 

\v2

{\sc {\bf 3.1.2.}  Two lemmas based on the  Cameron Martin formula.} \, It is consistent to our notation to write
\beqn\label{1dim}
q_a^{(1)}(x,t)=\frac{ P_x^{BM}[ T_a\in dt]}{dt} = \frac{x-a}{\sqrt{2\pi t^3}} e^{- (x-a)^2/2t} \quad (x>a).
\eeqn
Recall that  $q^\nu_a = q^{(2\nu+2)}_a$, and    $P_x^{BS(\nu)}$,  $P^{BM}_x$, $X_t$ and $T_a$   are introduced at the beginning of this section.
%
\begin{Lem} \label{lem2.1} \,  Put $ 
\b_\nu = \frac18 (1-4\nu^2)$ $(\nu\geq 0)$.  Then
\beqn\label{Eq0}
q_a^{\nu}(x,t) =q_a^{(1)}(x,t) \bigg(\frac{a}{x}\bigg)^{\nu+\frac12} E^{BM}_x\Bigg[\exp\bigg\{ \b_\nu \int_0^{t} \frac{ds}{X_s^{2}}\bigg\}\,\Bigg| T_a =t \Bigg].
\eeqn
\end{Lem}
\v2\n
\pf \, We apply the  formula of drift transform (based on the  Cameron Martin formula). Put    $Z(t) = e^{\int_0^t \ga(X_s)d X_s -\frac12 \int_0^t \ga^2(X_s)ds}$, where $\ga(x)=(\nu+\frac12)x^{-1}$  and  $X_t$ is a  linear Brownian motion.
Then 
\beqn\label{4.1.1}
\int_{t-h}^t q_a^{\nu}(x,s)ds= P^{BS(\nu)}_x[t-h\leq T_a <t] =E_x^{BM}[Z(t); t-h\leq T_a <t]
\eeqn
for $0<h<t$. 
By Ito's formula we have 
$\int_0^t dX_s /{X_s}= \lg (X_t/X_0) + \frac12\int_0^{t} ds/X_s^2$ ($t < T_0$).
Hence 
$$Z(T_a)= \bigg(\frac{a}{X_0}\bigg)^{\nu+\frac12} \exp\bigg[ \frac{1-4\nu^2}{8}\int_0^{T_a}\frac{ds}{X_s^2}\bigg],
$$
which together with (\ref{4.1.1}) leads to the identity (\ref{Eq0}). 

%
\begin{Lem} \label{lem2.2} \, For  $\la \geq 0$
\beqn \label{Eq1}
E_x^{BS(\nu)}\Bigg[\exp\bigg\{ -\frac{\la(\la+2\nu)}2  \int_0^{t} \frac{ds}{X_s^{2}}\bigg\}\,\Bigg| T_a =t \Bigg]
=\bigg(\frac{x}{a}\bigg)^\la\frac{q_a^{\la+\nu}(x,t)}{q_a^{\nu}(x,t)}. 
\eeqn
\end{Lem}
\v2\n
\pf\,  Write $\tau = \int_0^t X_s^{-2}ds$.  By the same drift transformation as applied in the preceding proof we see
$$\frac{E_x^{BS(\nu)}[e^{-\frac12 \la(\la+2\nu) \tau}; T_a \in dt]}{dt} 
= q_a^{(1)}(x,t) \bigg(\frac{a}{x}\bigg)^{\nu+\frac12} E^{BM}_x[e^{-\frac12 \la(\la+2\nu) \tau} e^{\b_\nu\tau}\,|\,T_a =t].$$
Noting $-\frac12 \la(\la+2\nu) +\b_\nu = \b_{\la+\nu}$  we apply  (\ref{Eq0}) with $\la+\nu$ in place of $\nu$ to see that the right-hand side above is  equal to
$(x/a)^\la q_a^{\la+\nu}(x,t)$, while the left-hand side is  equal to that of (\ref{Eq1}) multiplied  by $q^\nu_a(x,t)$, hence we have (\ref{Eq1}).  \qed

 \v2
 {\bf 3.1.3.}  {\sc Proof of Proposition  \ref{lem3.20}.}  For the proof we apply  the skew product representation of two-dimensional Brownian motion.  Let $Y(\cdot)$ be a standard linear Brownian motion with $Y(0)=0$ independent of $|B_\cdot|$. Then $\arg B_t  -\arg B_0$ has the same law as $Y(\int_0^t |B_s|^{-2}ds)$ (\cite{IM}), so that
 $$
E_{x{\bf e}}[ e^{i\la\arg B(\sigma_a)};\sigma_a \in dt]
 = E_{x{\bf e}}\otimes E^Y\Big[ e^{i\la Y\big({ \int_0^t |B_s|^{-2}ds}\big)};\sigma_a \in dt\Big]
 $$
 where $E^Y$ denotes the expectation with respect to the probability measure of $Y(\cdot)$ and $\otimes$ signifies the direct product of measures (with an abuse of notation). Note  that  
 $|B_t|$ is a two-dimensional Bessel process (of order $\nu=0$) and
   take    the conditional expectation of $e^{i\la Y\big( \int_0^t |B_s|^{-2}ds\big)}$ given $|B_{s}|, s\geq t$ to  find  the equality
 \beq
E_{x{\bf e}}[ e^{i\la\arg B(\sigma_a)} \,|\, \sigma_a = t] = E_x^{BS(0)} \bigg[\exp\Big\{-\frac{\la^2}2  \int_0^t X_s^{-2}ds\Big\}\,\Big|\, T_a =t\bigg],
\eeq
of which,  by  formula (\ref{Eq1}), the right-hand side  equals
$$(x/a)^{|\la|}q_a^{|\la|}(x,t)/q_a^{0}(x,t),$$
showing the required identity. \qed  

\v2

Let  $b>a$. Then  for each  $s>0$, the ratio $q^{(2)}_b (x,t-s)/q^{(2)}_a(x,t)$ is asymptotic to $\sqrt{b/a}\, e^{(b-a)v}e^{- \frac12 v^2 s}$ as $x/t\to v$, $t\to\infty$  and, on  considering the hitting of $U(b)$,  we  observe
\begin{eqnarray}
&& f_a(\th;v) d\th \nonumber \\
&&= \lim \frac{ \int_0^t \int_\R P_\x[\arg B_{\sigma_b} \in d\th'\,|\, \sigma_b =t-s] q^{(2)}_b(x,t-s) P_{be^{i\th'}} [\arg B_{\sigma_a} \in d\th, \sigma_a\in ds] }{q^{(2)}_a(x,t)} \nonumber\\
&&= \sqrt{\frac{b}{a}} \, e^{(b-a)v}\int_\R E_{b e^{i\th'}}\Big[e^{-v^2\sigma_a/2}; \arg B_{\sigma_a} \in d\th\Big]f_b(\th';v)d\th'
\label{f_v}
\end{eqnarray}
(with an appropriate interpretation of  $ \arg B_{\sigma_a}$ under $P_{be^{i\th'}} $),
which shows that $f_a(\th;v)>0$ for all  $\th$  and all $ v>0$.

\v2\n
{{\bf 3.2.} \sc An upper  bound of $q_1^{(2\la+2)}(x,t)$ for large $\la$.}  

For the proofs of Theorems \ref{thm1.1} and \ref{thm1.21} we need  a prtinent upper   bound of the characteristic function appearing  in Proposition \ref{lem3.20} for large integral values of  $\la $.  To this end we prove Lemma \ref{lem3.40}  below.  The result is extended to non-integral values of $\la$  in Lemma \ref{lem3.51} that verifies   the uniform convergence of the limit  appearing    in (\ref{f})  of   the conditional density for  $\arg B(\sigma_a)$.  
\v2
{\bf 3.2.1.} Here we prove the following  lemma.  
\begin{Lem}\label{lem3.40}   There exists  constants $C_1$ and $A_1>0$ such that for all  $n=1, 2, \ldots $, $t>1$ and  $x> 1$,
\beqn\label{q/p}
q_1^{(2n+2)}(x,t) \leq C_1 (A_1/n)^n p_{t+1/n}^{(2n+2)}(x).
\eeqn
\end{Lem}
\v2\n
\pf~
By the identity
$$p_{t+\e}^{(2n+2)}(x) = \int_0^{t+\e} q_1^{(2 n+2)}(x,t+\e -s) p^{(2n+2)}_{s}(1)ds$$
we have 
$$p_{t+\e}^{(2n+2)}(x) \geq  \bigg[\inf_{0\leq s\leq \e}  q_1^{(2 n+2)}(x,t + s)\bigg] \int_0^\e \frac{e^{-1/2s}}{(2\pi s)^{n+1}}ds$$
for every $0<\e<t$. We choose $\e=1/n$ and evaluate  the last integral  from below to see 
$$ \int_0^{1/n} \frac{e^{-1/2s}}{(2\pi s)^{n+1}}ds = \int_{n/2}^\infty e^{-u}u^{n-1}\frac{du}{2\pi^{n+1}}  \geq   \frac{A_0}{\sqrt n}\bigg(\frac{n}{e\pi}\bigg)^{n} $$
for some universal constant $A_0>0$. If $x>2$, we apply   the  inequality of Harnack type given in the next lemma to find  the inequality (\ref{q/p}). 

It remains to deal with the case  $1<x < 2$, which however  can be reduced to the case $x=2$. Indeed, by partitioning the whole event according as  $2$ is reached before  $t/2$ or not, we see (by recalling  $ q_1^{(2n+2)} = q_1^{n}$) that if $1<x<2$,
$$q_1^{n}(x,t) = P^{BS(n)}_x[T_1\wedge T_2>t/2]\sup_{1<y<2} q^{n}(y, t/2)  + \sup_{t/2 \leq  s < t} q_1^{n}(2,s).$$
The required upper bound of  the second term  on the right-hand side   follows from the result for $x=2$ since $p_{s}^{(2n+2)}(2) \leq 4^{n+1} p_{t+1/n}^{(2n+2)} (x)$ for $t/2 \leq  s<t$.
As for  the first term,   by Lemma \ref{lem2.1} we infer  that the  supremum involved in it  is bounded by a universal constant (since $\b_\nu\leq 0$ for $\nu\geq 1$). 
On the other hand, by the same drift transform that is  used in the proof of Lemma \ref{lem2.1} we see 
\beq
P^{BS(n)}_x\Big[T_1\wedge T_2> \frac{t}2\Big] &=& E^{BM}_x\Bigg[\bigg(\frac{X_{t/2}}{x}\bigg)^{n+\frac12} \exp\bigg\{ \b_n\int_0^{t/2} \frac{ds}{X_s^{2}}\bigg\}; T_1\wedge T_2>\frac{t}2\Bigg]\\
&\leq& e^{1/8}2^{n+1/2}e^{- n^2t/16} P^{BM}_x[T_1\wedge T_2> t/2],
\eeq
 which  is enough for the required bound. \qed

\begin{Lem}\label{lem-2}\,  There exist   constants $C_2>1$ and $A_2>0$  such that whenever  $x\geq 2$ and  $n=2, 3,\ldots$, 
$$q_1^{(n)}(x,t-\tau)\leq  C_2 A_2^n  q_1^{( n)}(x,t) \quad \mbox{for \,\, $t>1$ \, and  \, $0\leq \tau\leq 1/n$},$$
or, equivalently,\,\,  $q_1^{(n)}(x,t)\leq  C_2 A_2^n \inf_{ 0\leq s \leq 1/n}q_1^{( n)}(x,t+s) $ \, for \,  $t>1-1/n$. 
\end{Lem}
\v2\n
\pf \,  Let $Q$ be  the hyper-cube in $\R^n$ of side length $2$ and centered at the origin and put $D=\{(\y,s): \y\in Q, 0<s< 1+\tau\}$, the cubic cylinder with the base $Q\times \{0\}$ and of height $1+\tau$. The function $u(\y,s):=q_1^{(n)}(|\x+\y|,t-s)$ satisfies the  equation $\partial_s u + \frac12 \sum_{j=1}^n \partial_{j}^2 u=0$ in $D$, where $\partial_j$ denotes the partial derivative w.r.t. the $j$-th coordinate of $\y$.  Let  $p^0_s(x,y)$ be the heat kernel on the physical space  $[-1,1]$ with zero Dirichlet boundary and put 
$$p^0_s(\x,\y)=\Pi_{j=1}^n p_s^0(x_j,y_j) \quad  \mbox{and} \quad K({\bf S},s)=\pm \partial_ j p_s^0({\bf 0},\y)|_{\y={\bf S}},$$
where the sign is chosen so that $\pm \partial_ j $ becomes inner normal derivative at ${\bf S}\in \partial Q$. Then
$$u({\bf 0},\tau)= \int_{\partial Q}d {\bf S}\int_\tau^{1} K({\bf S}, s-\tau) u({\bf S},s)ds+ \int_Q p^0_{1-\tau}({\bf 0},\y)u(\y,1)d\y.$$
Since all the functions involved in these two  integrals are non-negative, we have 
$$q_1^{( n)}(x,t)= u({\bf 0},0) \geq \int_{\partial Q}d {\bf S}\int_\tau^{1} K({\bf S}, s) u({\bf S},s)ds+ \int_Q p^0_{1}({\bf 0},\y)u(\y,1)d\y,$$
and, comparing the right-hand side with the integral representation of $u({\bf 0},\tau)=q_1^{( n)}(x,t-\tau)$, we have
 $q^{( n)}(x,t-\tau) \leq M_n q^{( n)}(x,t)$, where $M_n = M'_{n}\vee M''_{n}$ with
 $$  M'_{n}=\sup_{{\bf S}}  \sup_{ \tau<s<1}  \frac{K({\bf S}, s-\tau) }{K({\bf S}, s) }, \quad M''_n= \sup_{\y}\frac{ p^0_{1-\tau}({\bf 0},\y)}{ p^0_{1}({\bf 0},\y)}.$$
 
 We must find a positive  constant $A_2$ for which   $M_n < C_2A_2^n$ if $\tau <1/n$. By the reflection principle we have
 $$p^0_s(0,y) =\sum_{k=-\infty}^\infty (-1)^k p_s^{(1)}(y-2k).$$
 Since $\sup_y  p^0_{1-\tau}(0,y) /p^0_1(0,y) $ tends to unity as $\tau \to0$,  we have
 $M''_n<2^{n}$ for all $\tau$ small enough. To find an upper bound of  $M_n'$ we
  deduce   the following bounds: for some  constant $C\geq 1$, 
 \beqn\label{Ineq0}
 \frac{p^0_{s-\tau}(0,y)}{p^0_s(0,y)}\leq C \sqrt{\frac{s}{s-\tau} }\quad \mbox{for} \quad \tau <s\leq 1,\,  |y| <1; \mbox{and}
 \eeqn
   \beqn\label{Ineq1}
   \frac2{s}p^{(1)}_s(1) - \frac6{s}p^{(1)}_s(3) < \mp \frac{\partial}{\partial y}p^0_s(0,y)\Big|_{y=\pm1} < \frac2{s}p^{(1)}_s(1)
  \quad \mbox{for} \quad 0 <s\leq 1.
  \eeqn
  The inequalities in  (\ref{Ineq1}) are  easy to show and its proof is omitted.   As for (\ref{Ineq0}) we  observe that if $\tau \leq s/2$, then   $s>s-\tau >s/2$ so that  the ratio on the left is bounded, while uniformly for   $|y|>1/2$   and for $\tau>s/2$, the ratio tends to zero   as $s\to 0$; in the remaining case   $|y|\leq 1/2$,  $ s/2< \tau<s$   the inequality (\ref{Ineq0}) is obvious. 
From  (\ref{Ineq0}) and (\ref{Ineq1}) we see that  if $s\geq  2\tau$, $M'_n <  (C 2^{1/2})^n$; also if  $\tau < s<2\tau$, then
 for $\tau =1/n$ small enough, 
 $$M'_n < 2\bigg[C \sqrt{\frac{s}{s-\tau} }\, \bigg]^{n}\exp\Big\{-\frac{\tau}{2(s-\tau)s}\Big\}\leq 2C^n\exp\Big\{-\frac{n}{2(u-1)u} + \frac{n}2  \lg\frac1{u-1}\Big\},$$
 where we put $u= s/\tau$. Thus,  putting  $m= \sup_{1\leq u\leq 2} \Big[- \frac{1}{2(u-1)u}+ \frac12 \lg \frac1{u-1}\Big]$ we have  $M'_n \leq 2 (Ce^m)^n$. The proof of the lemma is complete.
 \qed

\v2
{\bf 3.2.2.} {\sc  Convergence of the density for $\arg B_t$ conditioned on  $\sigma_a=t$ ($d=2$). }\, 
Here we  prove that the convergence in  (\ref{f}) holds uniformly in $\th$ locally uniformly in $v$. 
\begin{Thm}\label{thm3.3}  Let $d=2$. For each  $M>1$, uniformly for 
 $\th\in \R$ and $x \in (a, Mt)$, as $t\to\infty$
 $$\frac{P_{x\be}[\arg B(\sigma_a)\in d\th \,|\, \sigma_a =t]}{ \,d\th}= f_a(\th;  x/t)(1+o(1)).$$
 \end{Thm}
 
  For the proof   we need  the following extension of Lemma \ref{lem3.40}.
 \begin{Lem}\label{lem3.51} There exist   constants $C$ and $A>0$ such that for all  $\la > 1$, $t>1$ and $x> 1$,
\beqn\label{q/p2}
q_1^{(2\la+2)}(x,t) \leq C(t/x)^\de  (A/\la)^\la p_{t+1/\la}^{(2\la+2)}(x).
\eeqn
Here  $\de$ denotes the fractional part of $\la$. 
\end{Lem}
\v2\n
\pf~  We may and do suppose $\la\in (n, n+1)$ for a positive integer  $n$.  Remember that   $(P_x^{BS(\nu)}, X_t)$ designates a Bessel process  of dimension $2\nu+2$.  Put $\de =\la-n$ and  $\ga(y) = \de/y$. Then the drift transform gives 
$$P_x^{BS(\la)}[ \Ga; T_1 \geq t] =  E_x^{BS(n)} [Z(t); \Ga, T_1 \geq t]$$
 for any event $\Ga$ measurable w.r.t. $(X_s)_{s\leq t}$ (cf. e.g. \cite{IW}). Here, since the drift term of $X_t$ under $P^{BS(n)}$ equals $(2n+1)/2X_t$,  
$$Z(t)  =  \int_0^t \ga(X_s)dX_s - \frac12 \int_0^t\Big[ \frac{\ga(X_t)(2n+1)}{X_t} + [\ga(X_s)]^2\Big]ds.$$
 By Ito's formula
we have $\int_0^t \ga(X_s)dX_s = \de\lg (X_t/X_0) - \frac12\int_0^t \ga\,'(X_s)ds$. Observing  $\ga\,'(y)+(2n+1)\ga(y)/y +\ga^2(y) = (2n\de +\de^2)/y^2$, as in Section {\bf 3.1.2} we find 
$$q_1^{(2\la+2)}(x,t) = x^{-\de} E_x^{BS(n)}[e^{-\frac12 (2n\de+\de^2)\int_0^t X_s^{-2}ds}\, |\,  T_1=t]q_1^{(2n+2)}(x,t). $$  
The conditional expectation being dominated by  unity,   substitution from Lemma \ref{lem3.40} yields
$$q_1^{(2\la+2)}(x,t)  \leq (t/x)^\de t^{-\de} C_1 (A_1/n)^n  p_{t+1/n}^{(2n+2)}(x)\leq C (t/x)^\de n^\de (A_1/\la)^\la 
p_{t+1/\la}^{(2\la+2)}(x),$$
showing  the inequality of the lemma  with  any $A>A_1$. \qed

\v2\n
{\it Proof of Theorem \ref{thm3.3}.}
Let $a=1$. From the lemma above and Proposition \ref{lem3.20} we see that for $t>1, x>1$, 
 $$ E_{x{\bf e}}[ e^{i\la\arg B(\sigma_a)}|\sigma_1 = t] 
\leq C' \Big[1\vee \lg\frac{t}{x}\Big]^{2}\Big(\frac{A}{|\la|}\Big)^{|\la|} \bigg(\frac{x}{t}\bigg)^{|\la|-\de} \exp \frac{(x/t)^2}{2|\la|},  \,\, 
\quad n< |\la| \leq n+1, $$
where we have  also used the lower bound $q^{(2)}_1(x,t)\geq C[1\vee \lg(t/x)]^{-2}\, p_t^{(2)}(x)$.  On recalling  (\ref{K/nu}) as well as  Theorem \ref{prop2}  this  
shows that  the characteristic function on the left converges to  $K_0(x/t)/K_\la(x/t)$, the Fourier transform of  $f_1(\cdot, x/t)$,  in $L_1(d\la)$  uniformly for $x/t <M$, hence the uniform  convergence asserted in  the lemma. \qed

\v2\n
   {\sc {\bf 3.3.}     Distribution of  $\Th(\sigma_a)$}.\, 
   
In this subsection we give proofs of  Theorems \ref{thm1.1} and \ref{thm1.21}.  To facilitate  the exposition we first introduce the  conditional density $g(\th;x,t) $. We  then  expand $g$  into  a series of spherical functions which  almost immediately leads to  (a refined version of)  Theorem \ref{thm1.1}  and  to Theorem \ref{thm1.21}.
   \v2
    {\sc {\bf 3.3.1.} The  conditional density $g(\th;x,t) $.} Let $\th(\xi)$ denote the colatitude of  $\xi\in \partial U(a)$ as before. 
   By rotational symmetry  around the axis  $\eta\be, \eta\in \R$  we can define  $g(\th;x,t)$ by 
\beqn\label{sin0}
g(\th;x,t) := \frac{P_{x\be}[B(\sigma_a)\in d\xi \,|\, \sigma_a =t] }{m_a(d\xi )}, \quad  \th = \th(\xi) \in [0,\pi].
\eeqn
Denote  the colatitude  $\th(B_t)$ by  $\Th_t \in [0,\pi]$,  
 so that
 $$ \cos \Th_t = {\bf e}\cdot B_t /|B_t|.$$ 
 
 Let  $d\geq 3$ and  $d\xi  = a^{d-1}\sin^{d-2} \th\, d o\times d\th$, where  $d o$ designates a $(d-2)$-dimensional surface element of $(d-2)$-dimensional unit sphere. 
 Then  $m_a(d\xi) = \sin^{d-2} \th \,d\th |do|/\omega_{d-1}$ and we see that 
\beqn\label{sin}
g(\th;x,t)=  \frac{P_{x\be}[\Th(\sigma_a)\in d\th \,|\, \sigma_a =t]}{ \mu_d^{-1}\sin^{d-2}\th \,d\th}.
\eeqn
Here $\mu_d= \int_0^\pi \sin^{d-2}\th d\th = \omega_{d-1}/\omega_{d-2}$. 

When $d=2$,  we have $\Th_t = |{\rm Arg} \,B_t|$ and 
$$g(|\th|;x,t) =  2\pi \frac{ P_{x\be}[{\rm Arg} \, B(\sigma_a)\in d\th \,|\, \sigma_a =t]}{ d\th}, \quad \th\in (-\pi,\pi).$$
Thus the measure $g(|\th|;x,t) d\th/2\pi$ on $|\th|\leq \pi$ is  the probability  law of ${\rm Arg}\, B(\sigma_a)$  under $P_{x\be}[\cdot\,|\, \sigma_a=t]$ and we may/should  naturally regard $g(|\th|;x,t)$ as a (continuous) function on the torus  $\R/2\pi \Z \cong [-\pi, \pi]$.
 It is noted  that by letting $\om_0=2$ so that $\mu_2=\pi$,  the last expression  conforms to  (\ref{sin}).
(In  (\ref{sin})  the differential quotient at the end point $\th =0$ (or $\pi$)     is understood to be  the right (resp. left)  derivative of the distribution function.)


\v2
{\bf 3.3.2.} {\sc  Series expansion of $g(\th; x,t)$  when $d=2$}.\, 
Let  $d=2$ and $g(\th, x,t)$ be given as above. Denote by
  $\a_n=\a_n(x,t)$, $n=0, 1, 2,\ldots$  the  coefficients of the  {\it Fourier cosine series} of  $g(\th) = g(\th;x,t)$, $\th\in [0, \pi]$: $\a_0=\pi^{-1}\int_0^\pi g(\th)d\th = 1$ and for $n\geq 1$, 
 $$\a_n = \frac2{\pi}\int_0^\pi g(\th)\cos n\th\, d\th = 
 2 E_{x\be}[\cos n \Th(\sigma_a)\,|\, \sigma_a =t],$$
so that 
\beqn\label{F-exp}
g(\th; x,t)= \sum_{n=0}^\infty \a_n(t,x)   \cos n\th,
\eeqn
where the Fourier series  is uniformly convergent (with any  $x, t$ fixed)    as one may infer from the smoothness of $g$ (or alternatively  from our estimation of $\a_n$ given in (\ref{alpha})    below). Since  $E_{x\be}[\cos n \Th(\sigma_a)\,|\, \sigma_a =t] = E_{x\be}[e^{in \arg B(\sigma_a)}\,|\, \sigma_a =t]$,
substitution from   Proposition  \ref{lem3.20}  yields
\beqn\label{-0}
\a_n(x,t) = 2\frac{q_a^{(2n+2)}(x,t)}{q_a^{(2)}(x,t)} \bigg(\frac{x}{a}\bigg)^{n}. 
 \eeqn
Based on  this formula  we derive the next  result that provides an exact asymptotic form of the error term in Theorem  \ref{thm1.1} (i).  
    (As another possibility  one may use  a classical formula for $g(\th;x,t)q^{(2)}_a(x,t)$  that we give   in Appendix.)   
\begin{Thm} \label{thm-1} \,  Let $d=2$. Uniformly for $\th\in [0, \pi]$ and $x>a$, as $t\to\infty$ with $x/t\to 0$,
$$g(\th; x,t) = 1 + \frac{ax}{t} \ell_0(x,t)\Big[(1+o(1))\cos \th + O\Big(\frac{x}{t}\Big)\Big],$$
where 
$$\ell_0(x,t)=\Big(1-\frac{a^2}{x^2}\Big) \frac{(\lg t)^{2}}{2 \lg (x/a)} \quad \mbox{ if }\,\, 1<x<\sqrt t ; \mbox{ and}  \,\, = 2 \lg \frac{t}{x} \quad  \mbox{  if} \,\,\,  x>\sqrt t.$$
\end{Thm}
\v2\n
\pf \,  By   elementary computation we deduce from   Theorem A and (\ref{-0}) that 
\beqn\label{-3}
\a_1= \frac{ax}{t} \ell_0(x,t)(1+o(1)) 
\eeqn 
 as $ x/t \to 0$. Plainly   $\a_n(x,t)\geq 0$.  It therefore  suffices to show that 
\beqn\label{-1}
\sum_{n=2}^\infty \a_n(x,t) = O\bigg(\frac{x^2}{t^2} \ell_0(x,t)\bigg).
\eeqn
Although Theorem A also  yields $ \a_n(x,t) = O\Big((x/t)^n\ell_0(x,t)\Big)$ for each $n= 2, 3 \ldots$,   for the present purpose  we need  an upper bound valid  uniformly  in $n$. Such a uniform bound    is provided by 
 Lemma \ref{lem3.40} and on using  it  
\beqn\label{alpha}
\a_n(x,t) \leq C_2 \frac{A_1^n}{n^n} \bigg(\frac{x}{t}\bigg)^{n} \frac{e^{-x^2/2t}}{tq_1^{(2)}(x,t)} \leq C_3 \frac{A_1^n}{n^n} \bigg(\frac{x}{t}\bigg)^{n}  \ell_0(x,t),
 \eeqn
which   implies  (\ref{-1}).   
    \qed

\v2
{\bf 3.3.3.} {\sc  Series expansion of $g(\th; x,t)$ when $d\geq 3$.}
\, Recall 
$$ P_{x\be}[\Th(\sigma_a)\in d\th\,|\, \sigma_a =t] =   \mu_d^{-1}g(\th;x,t)  \sin^{d-2} \th\,d\th.
$$
\begin{Thm} \label{thm-11} \,  Let $d\geq 3$. For $\th\in [0, \pi]$ and $x>a$,
$$g(\th;x,t) = \sum_{n=0}^\infty\bigg(\frac{x}{a}\bigg)^n\frac{q_a^{n+\nu}(x,t)}{q_a^{\nu}(x,t)} h_n(0)h_n(\th),$$
where  $h_n(\th)$ denotes the $n$-th normalized eigenfunction  of  the Legendre process of order $\nu$ (see Section 6). 
\end{Thm}
\v2\n
\pf \, Let $(P^{L(\nu)}_\th, \Th_t)$ denote the Legendre process (on the state space $[0,\pi]$) of order $\nu$.  Then by the skew product representation of $d$-dimensional Brownian motion we have
\beq
P_\x[\Th(\sigma_a) \in d\th, \sigma_a \in dt]
=(P_{\th_0}^{L(\nu)}\otimes P_x^{BS(\nu)})[\Th_\tau \in d\th\,|\, T_a =t]q^{(d)}(x,t),
\eeq
where $\tau = \int_0^{T_a} X_s^{-2}ds$ and $\th_0$ is the colatitude  of  $\x$. We apply  the spectral expansion of the density of the distribution of $\Th_t$ (see (\ref{spexp})) and Lemma \ref{lem2.2} in turn to deduce that
\begin{eqnarray}\label{hexp}
&& (P_{\th_0}^{L(\nu)}\otimes P_x^{BS(\nu)})[\Th_\tau  \in d\th\,|\, T_a =t] /d\th \nonumber\\
&&= E_x^{BS(\nu)}\Bigg[ \sum_{n=0}^\infty \exp\Big\{-\frac{n(n+2\nu)}{2} \tau \Big\}h_n(\th_0)h_n(\th) \frac{\sin^{d-2} \th}{\mu_d} \,\Bigg| \, T_a =t \Bigg] \nonumber\\
&&= \frac{1}{\mu_d}\sum_{n=0}^\infty\bigg(\frac{x}{a}\bigg)^n\frac{q_a^{n+\nu}(x,t)}{q_a^{\nu}(x,t)} h_n(\th_0)h_n(\th) \sin^{d-2} \th.
\end{eqnarray}
Comparing this  with  (\ref{sin})  shows the formula of the theorem.    \qed
\v2
In view of the defining  identity  (\ref{sin0})   the case $d\geq 3$ of Theorem \ref{thm1.1} follows from 
\begin{Cor} \label{cor-12} \,  Let $d\geq 3$. Uniformly for $\th\in [0, \pi]$ and $x>a$, as $x/t\to 0$
$$g(\th;x,t) = 1 + \frac{ax}{t}\bigg[\frac{1-(a/x)^{d} }{1-(a/x)^{d-2}}\Big(\frac{d}{d-2}+o(1)\Big)\cos \th + O\Big(\frac{x}{t}\Big)\bigg].$$
\end{Cor}
\v2\n\pf\, The asserted formula  is derived as in the case $d=2$ by observing that  $h_1(0)h_1(\th) = 2(\nu +1)\cos \th$ (see Section {\bf 6.1.1}) and 
$$\frac{x}{a}\frac{q_a^{1+\nu}(x,t)}{q_a^{\nu}(x,t)}  \sim \frac{ax}{t}\cdot \frac{1-(a/x)^{2+2\nu}}{ 2\nu(1-(a/x)^{2\nu})}(1+o(1)).$$
\qed

\v2
{\bf 3.3.4.} {\sc Proof of Theorem \ref{thm1.21}.} \,
Proof of Theorem  \ref{thm1.21} proceeds as follows.    For  $d=2$  Theorem \ref{thm-11}  is valid with  $h_n(0)h_n(\th)$ replaced by  $2\cos n\th$ if $n\geq 1$  as we  have already   observed  (see (\ref{F-exp}) and  (\ref{-0})); here it  is  warned that if $d=2$ the product $h_n(\th_0)h_n(\th)$ must be replaced by $2\cos n(\th-\th_0)$ ($n\geq1$)  in (\ref{hexp}).  In any case  substitution from  (\ref{K/K}) gives  the relation of Theorem \ref{thm1.21} for $d=2$ at a formal level.
The relation  (\ref{K/K})  is immediately extended to 
$$\frac{q_a^{|\la|+\nu}(x,t)}{q_a^\nu(x,t)} \bigg(\frac{x}{a}\bigg)^{|\la|}
\sim  \frac{K_{\nu}(av)}{K_{\nu+|\la|}(av)}\qquad (x/t\to v).  
$$
With these remarks  as well as  (\ref{q/p}) taken into account we obtain  from  (\ref{F-exp})  and Theorem \ref{thm-11} that for all $d\geq 2$,
as $x/t\to v$
$$g(\th;x,t)  - \sum_{n=0}^\infty \frac{K_\nu(av)}{K_{\nu+n}(av)} b_nh_n(\th) \, \longrightarrow\,  0, 
$$
uniformly in $\th\in [0,\pi]$ and $0<v<M$ for each $M$. Here $b_nh_n(\th)= 2\cos n\th$ for $n\geq 1$ if  $d=2$ and $b_n=h_n(0)$ if $d\geq 3$.
This shows  Theorem \ref{thm1.21}  except for identification of the constant factor  in the case $d\geq 3$,  which we give at the last line of Section {\bf 6.1.1}.

\v2\n
{\sc Remark 3.} \,  There exists an unbounded and increasing positive  function  $C(v)$, $v>0$ such that $C(0+)\geq 1$ and 
$$1/C(x/t) \leq g(\th;x,t) \leq C(x/t) \qquad  (0\leq \th\leq \pi, t>1).$$
The upper bound follows from   (\ref{F-exp}), 
Theorem \ref{thm-11} and estimates like (\ref{alpha}), while the lower bound can be verified by an argument analogous to the one as 
made at (\ref{f_v})  (or  in Section {\bf 5.4}).


\section{ Estimates of the hitting density for $t<1$ }
 Put for $z>a$
\beqn\label{h_z00}
h_a(z,t,\phi) = \frac{P_{z\be}[\Th(\sigma_a)\in d\phi,   \sigma_a \in dt]}{ \mu_d^{-1}\sin^{d-2}\phi\, \,d\phi  dt}, \quad \phi\in [0,\pi],
\eeqn
or,  by means of  $g=g_a$  given in (\ref{sin}), 
$$ h_a(z,t,\phi)= g_a(\phi; z, t) q^{(d)}_a(z,t);$$
recall that $g_a(\phi;z, t)$ represents the density with respect to $m_a(d\xi)dt$ evaluated at  $\xi$ with colatitude $\phi = \th(\xi)$  of the hitting  site distribution  conditional on $\sigma_a=t$ and $B_0=z\be$. In view of  rotational symmetry of Brownian motion it follows that
for any $\xi\in \partial U(a)$ with  $\z\cdot\xi /xa=\cos \th$ and $\z \notin U(a)$,  
\[
h_a(|\z|,t, \th) := \frac{P_{\z}[B(\sigma_a)\in d\xi,  \sigma_a \in dt] }{m_a(d\xi )dt}.
\]

In this section we provide some upper and lower bounds of $h_a(z,t,\phi)$  for $t<1$, which are used    in the next section for estimation of  it  when  $z/t$ along with  $t$ tends to infinity. We include  certain  easier results for $t\geq 1$. The main results of this section are   
 given in  Lemmas \ref{lem3.5} and  \ref{Imp1}.
 For all dimensions $d\geq 2$ the function $h_a(z,t,\phi)$ satisfies the scaling relation 
$$h_a( z, t, \phi) = a^{-2} h_1(z/a, t/a^2, \phi).$$

Throughout  this section $X_t$  always 
denotes a standard  linear Brownian motion. As in the preceding section    $P_y^{BM}$ and $E^{BM}_y$  denote   the  probability and expectation for  $X_t$,  and   $T_y$  the first passage time of $X$  to  $y$.  We shall apply  the skew product representation  of $d$-dimensional Brownian motion  and  the Bessel processes  of dimensions $d\geq 2$ will become relevant.  However,   most  of the results  of this section that  actually concerns   the Bessel processes   follows from  the   one for the linear Brownian motion $X_t$  because of  the boundedness of  the Radon-Nikodym density $Z(t)$ ($t<1$) that is given  in the proof of Lemma  \ref{lem2.1} (see Remark 4 below
for more details).
 
 \v2\n
{\bf 4.1.}  {\sc Some Basic Estimates.}
    \v2
\begin{Lem}\label{lem3.2i} \,   Let $b>0 $. For $0<y<b$  and $0< t \leq b^2$, 
$${\displaystyle \, \frac{P^{BM}_y[T_0\in dt, T_b< T_0]}{dt} \leq C\frac{yb^2}{t^2}p^{(1)}_t(b)}.$$
\end{Lem}
\v2\n
\pf \, By reflection principle it follows that 
$$\frac{P^{BM}_y[T_0\in dt, T_0< T_b]}{dt}  =  
\frac1{\sqrt{2\pi t^3}} \sum_{n= -\infty}^\infty( 2nb +y) \exp\bigg\{-\frac{ (2nb +y)^2}{2t}\bigg\}
$$
(\cite{KS}, (8.26)).  Writing the right-hand side above in terms of  $q_a^{(1)}$ (cf.  (\ref{1dim})) we see that
 \beq
\frac{P^{BM}_y[T_0\in dt, T_b< T_0]}{dt}  & =& q^{(1)}_0(y,t) - \frac{P^{BM}_y[T_0\in dt, T_0< T_b]}{dt} \\
&=&
 \sum_{n=1}^\infty [q^{(1)}_0(2nb-y,t)- q^{(1)}_0(2nb+y,t)].
 \eeq
On using the  mean value theorem   the difference under the summation symbol is dominated  by 
$$\frac{2y}{\sqrt{2\pi t^3}} \frac{[(2n+1)b]^2}{ t} e^{- [(2n-1)b]^2/2t} \quad (0<y < b, 0<t<b^2).$$
By easy domination of these  terms for $n\geq 2$ we find 
 the upper  bound of the lemma.
  \qed
 \v2\n
{\sc Remark 4.} \,Lemma \ref{lem3.2i} is extended to $d$-dimensional Bessel Processes  $|B_t|$  with essentially the same  bound if  the positions $0, y$ and $b$ are raised to  $a$, $a+y$ and $a+b$, respectively,    by using  the drift transformation. For later reference here   we give it in the form 
        \beqn\label{drift}
     P_{(a+y)\be}[ A \,|\, \sigma_{a} = t] = c_a(y,t) E_{a+y}^{BM}[ e^{\b_\nu \int_0^t X_s^{-2}ds}; A^X \,|\, T_a =t],
     \eeqn
    where $\b_\nu =\frac18 (1-4\nu^2)= \frac18(d-1)(3-d)$, $A$ is an event of the process $|B_s|, 0\leq s\leq t$,  $A^X$ the corresponding one for $X$ and
 $$c_a(y,t)  := \bigg(\frac{a}{a+y}\bigg)^{(d-1)/2}\frac{q_a^{(1)}(a+y,t) }{q_a^{(d)}(a+y,t)} =1 + O\bigg(\frac{t}{a(a+y)}\bigg) \quad (0<t< a^2,  y>0).$$
(The last equality follows from Theorem B.) 
\begin{Lem}\label{lem3.2ii} \,  For $\a >0$    there is  a  constant $\k_{\a,d}$ (depending  on  $d, \a$) such that  $$E_{(1+y)\be}\bigg[\bigg(\int_0^t\frac{ds}{ |B_s|^{2}}\bigg)^{-\a} \,\bigg|\,  \tau_{U(1+\la)}<t, \sigma_1=t\bigg] \leq \k_{\a,d} (1+\la^{2\a}) t^{- \a}$$ 
for $\la>0$, $0<y < \la$ and $  t<\la^2$,
 where $\tau_{U(b)}$  denotes the first exit time from $U(b)$.
\end{Lem}
\v2\n
\pf\, The proof is given only for  the case $d=1$.
Put $M_t = \max_{s\leq t} X_s$. Then  the conditional expectation in the lemma multiplied by  $t^\a$  is at most
   $$E^{BM}_{y} [ (1+M_t)^{2\a}\,|\, T_\la< T_0=t] \leq 4^\a  + 4^\a\frac{E^{BM}_y[M_t^{2\a}; T_\la<t \,|\, T_0=t]}{P^{BM}_y[ T_\la <t\,|\, T_0=t]}.$$
  The last ratio may be expressed as a weighted  average of  $ E^{BM}_\la[ M_{t-s}^{2\a}\,|\, T_0 =t-s]$ over $0\leq  s\leq t$, which, by virtue of scaling property, is dominated by   $C'_\a\la^{2\a}$, yielding the desired bound.  \qed
\v2

  \begin{Lem}\label{lem3.4} \,  There exists a constant $\k_{d}$ depending only on $d$ 
    such that  for  $0 < \la \leq 8$,
  $$h_a(a+y,t, \phi) \leq \k_{d} \frac{a^{2\nu+1}y}{t}  \bigg(p_t^{(1)}(y) p_t^{(d-1)}(a\phi) + \frac{ (\la a)^2}{t}p_t^{(d)}(\la a)\bigg)$$
whenever  $0\leq \phi<\pi$, $0< y< \la a$ and $ 0<t<(\la a)^2$.
  \end{Lem}
  \v2\n
  \pf\, We may let $a=1$.   Suppose $d=2$.  Let   $(P^Y, (Y_t))$  be a  standard Brownian motion  on the torus  
   $\R/ 2\pi \Z$ (identified with $(-\pi,\pi]$)  that is started  at 0 and  independent of $(B_t)_{t\geq 0} $. Then by skew product representation of $B_t$
    \beqn\label{h_skew}
h_1(1+y,t, \phi)= 2\pi (P^Y\otimes P_{(1+y)\be}) [ Y_\tau\in d\phi, \sigma_1\in dt ]/d\phi dt,
   \eeqn
   where $ \tau = \int_0^{\sigma_1} |B_s|^{-2} ds$  and   $\otimes$ signifies the direct product of measures. 
  We rewrite this identity by means of the linear Brownian motion $X_t$ only. Because of translation invariance of the law of  the increment of  $X_t$ we shift the starting point of $X_t$ so that $X_0=y$ and define
   \beqn\label{tau_1}
   \tau^X = \int_0^{T_0}\frac{ds}{(1+X_s)^{2}}.
   \eeqn
 We perform the integration of $Y$ first  and   apply  the drift transform (as in the proof of Lemma \ref{lem2.1})   to deduce from (\ref{h_skew}) that 
  \beqn\label{h_skew2}
h_1(1+y,t, \phi)= \frac{2\pi}{\sqrt{1+y}} E^{BM}_{y} \Big[e^{\frac18 \tau^X } p_{\tau^X}^{\rm trs}(\phi) \,\Big|\, T_0=t\Big] q_0^{(1)}(y,t),
   \eeqn
   where  $p_t^{\rm trs}(\phi)$ denotes the density of the distribution of $Y_t$.
We  break the conditional expectation  above  into two parts according as $T_0 < T_\la$ or  $T_0> T_\la$,  and  denote  the corresponding ones by $J(T_0<T_\la) $ and $J(T_0> T_\la )$, respectively.  Note that $\tau^X <t $ (under  $T_0=t$)  so that $ p_{\tau^X}^{\rm trs}(\phi) \leq C p_{\tau^X}^{(1)}(\phi)$ if  $\sqrt t  < \la\leq 8$. 
 Then,  using Lemma \ref{lem3.2i}  (with $b=\la$) and Lemma \ref{lem3.2ii} (with $\a= 1/2$) we observe
 \begin{eqnarray}\label{J1}
  J(T_0> T_\la ) 
 &=&  E_y^{BM} [e^{\frac18 \tau^X} p^{{\rm trs}}_{\tau^X}(\phi) \,|\, T_0=t> T_\la]  \times P^{BM}_{y}[ T_0>T_\la \,|\, T_0=t] 
 \nonumber\\
 &\leq& C e^{\frac18 t}E_y^{BM} [(\tau^{X})^{-1/2} \,|\, T_0=t> T_\la]  \times P^{BM}_{y}[ T_0>T_\la \,|\, T_0=t] 
\nonumber\\
&\leq& Ce^{\la^2/8}\frac{(1+\la)}{t^{1/2}}\times \frac{y\la^2}{t^2}p^{(1)}_t(\la) \times \frac{1}{q_0^{(1)}(y,t)}.
 \end{eqnarray}
 On the other hand, the trivial domination  
$P^{BM}_{y}[ T_0<T_\la \,|\, T_0=t] \leq 1$  yields
 \begin{eqnarray}\label{J2}
  J(T_0<T_\la) 
&\leq & Ce^{\la^2/8}E^{BM}_{y}[p_{\tau^X}^{(1)}(\phi) \,|\, T_0=t< T_\la] \\
   & \leq&C  e^{\la^2/8}(1+\la) p^{(1)}_t(\phi).  \nonumber
\end{eqnarray}
Here the second inequality  is due to the inequality  $p^{(1)}_{\tau^X}(\phi) \leq (1+\la)p^{(1)}_t(\phi)$ that is valid if  $(1+\la)^{-2}t<\tau^X<t$, hence   if  $t< T_\la$. On recalling $q_0^{(1)}(y,t) = (y/t)p_t^{(1)}(y)$
these together show the estimate of the lemma when $d=2$.  

The higher-dimensional case $d\geq 3$ can be dealt with  in the same way  in view of   what is  noted in Remark 4   and the fact that   the transition density of a (spherical) Brownian motion  on the $(d-1)$-dimensional sphere  is 
comparable with that on the flat space if $t$ is small (cf.  Section  {\bf 6.1.2}). The details are omitted.  \qed
  \v2
  \begin{Lem}\label{lem3.41} \,  Uniformly for $y>0$,  as  $(y^3+|\phi|^3 ) /t \to 0$ and  $t\downarrow 0$
  $$\frac{h_a(a+y,t, \phi)}{2\pi} = \frac{a^{2\nu+1}y}{t} p_t^{(1)}(y) p_t^{(d-1)}(a\phi)(1+o(1)).$$
  \end{Lem}
  \v2\n
  \pf\, This proof is performed by examining the preceding one. We suppose  $d=2$ and $a=1$. By virtue  of the identity (\ref{h_skew2}) it suffices to show
   \beqn\label{EXPT}
    E^{BM}_{y} \Big[e^{\frac18 \tau^X } p_{\tau^X}^{\rm trs}(\phi) \,\Big|\, T_0=t\Big]  = p_t^{(1)}(\phi)(1+o(1))
    \eeqn
 in the same limit as   in the lemma.  Given  $t>0$ we  put  $\la= \la(t)= t^{1/3}$. With $b= \la(t)$ the inequality of Lemma \ref{lem3.2i} holds true, hence   also  (\ref{J1}) and (\ref{J2}) do even though $\la(t)$ depends on $t$. From the constraint on $\phi, y$ and $t$ imposed in the lemma  it follows that
 \beqn\label{ratio}
 \frac{y+|\phi| +\sqrt t}{\la(t)} \to 0\quad \mbox{and}\quad \frac{\phi^2\la(t)}{t} \to 0.
 \eeqn
As before we break the expectation into two parts. The part $J(T_0> T_\la )$  is negligible, for  the last member in (\ref{J1}) is at most a positive multiple of  $t^{-3/2}p_t^{(1)}(\la)/p_t^{(1)}(y)$ and the latter is  $o(p_t^{(1)}(\phi))$ under (\ref{ratio}). As for 
  $J(T_0< T_\la )$   the  estimate from above  is provided by (\ref{J2}). For, $C$ in (\ref{J2})  that comes in  from the bound   $ p_{\tau^X}^{\rm trs}(\phi) \leq C p_{\tau^X}^{(1)}(\phi)$  may be taken arbitrarily close to 1 as $\tau^X<t \to 0$. The  estimate from below is obtained  by observing that if $T_0< T_\la$ (so that $\tau^X > (1+\la)^2t$), then 
 $$\frac{p_{\tau^X}^{(1)}(\phi)}{p_t^{(1)}(\phi)} = \sqrt{\frac{t}{\tau^X}} \exp\Big\{-\frac{\phi^2}{2t\tau^X}\int_0^t\frac{2X_s+X_s^2}{(1+X_s)^2}ds\Big\}\geq \frac{1}{1+ \la}e^{- 2\phi^2\la/t} \to 1.$$ 
The proof of the lemma is complete.  
  \qed

\v2
The estimate of Lemma \ref{lem3.4}, which  concerns the case when $(z-a)/t$ is bounded above, will be improved in  Lemma \ref{Imp1} of the next subsection. The next lemma   
  provides  a bound of $ h_a(z,t, \phi)$  valid  for a wide range of the variables $z$, $\phi$ and $t$.  To simplify  the description of it as well as of its proof  we bring in a notation  that represents $ h_a(z,t, \phi)$ in a different way.

For $\z\notin U(a)$, put
\beqn\label{h_z0}
h^*_a(\z,t) = \frac{P_{\z}[ B(\sigma_a)\in d\xi,   \sigma_a \in dt]}{ m_a(d\xi)  dt}\bigg|_{\xi=a\be},
 \eeqn
which   may be also understood to be  the  density evaluated at $(0,t)$ of the joint law of $(\Theta(\sigma_a), \sigma_a)$   under $P_\z$. 
  If  $\z\cdot\be/z=\cos \phi\neq -1$, $z=|\z|$, then  $ h_a(z,t, \phi)= h^*_a(\z,t)$ due to rotational symmetry of Brownian motion. 
   When  $d=2$ 
 these may be given as follows: 
 \beqn\label{h_z2}
h^*_a(ze^{i\phi},t) =   h_a(z, t, \phi) = 2\pi \frac{P_{z\be}[{\rm Arg}\, B(\sigma_a)\in d\phi,   \sigma_a \in dt]} {d\phi  dt}.
 \eeqn

  \v2
  \begin{Lem}\label{lem3.5} \, Let  $|\z|>a$ ($\z\in \R^d$)  and put $r=|\z- a\be|$.  Then for some  constant $\k_d$,
 $$\quad  h^*_a(\z,t) \leq \k_d q^{(d)}_a(z,t) \qquad \mbox{ if} \quad  t> a^2\vee ar; \,\mbox{and}$$
     $$h^*_a(\z,t) \leq \k_d a^{2\nu} \, \frac{a r}{t}\, p_t^{(d)}(r) \quad\, \mbox{if}   \quad  t \leq a^2\vee ar. $$
       \end{Lem}
  \v2\n
  \pf\,    The case  $ t \geq ar$  is readily disposed of.   Indeed    
  the asserted inequality is implied by   Theorems \ref{thm1.1} and \ref{thm1.21} (in conjunction with  Theorem A)  if $ t \geq a^2\vee ar$, and 
 by Lemma \ref{lem3.4} if $ar<t <a^2$ (note that   $p^{(d)}_t(r)\asymp p^{(d)}_t(0)$ in the latter case).

 In the rest of  proof we let $a=1$   and suppose $t\leq r$, the case  which plainly entails  $t<1\vee r$ and thus concerns the second bound of the lemma.   Take positive numbers  $\e<1$ and $R$  so that  $r-\e>R>\e$. Then, on considering the ball  about  $(1-\e)\be$ of radius $\e$,  
  \begin{eqnarray}\label{951}
  h^*_1(\z,t) &\leq& \e^{-2\nu-1}  h^*_\e( \z- (1-\e)\be, t)  \nonumber\\
   & =&\int_0^t  \frac{h^*_\e(\xi,t-s)}{\e^{2\nu+1}} \int_{\partial U(R)}  P_{\z- (1-\e)\be}[\sigma_{U(R)} \in ds, B_{s} \in d\xi]. 
  \end{eqnarray}
Here, in the middle member  we have the factor $\e^{-2\nu-1}$ in front of  $h^*_\e$  since  the uniform probability  measure of the surface element $d\xi$ at $\be$ of  the sphere $\partial U(\e)+ (1-\e)\be$   equals   $\e^{-2\nu-1}m_1(d\xi')$ with  $d\xi' \subset \partial U(1)$ designating  the projection of $d\xi$ on  $\partial U(1)$ (see Remark 5 following this proof for the inequality).
    Write 
  $$r_*=r_*(\e) = | \z -(1-\e)\be|, \,\,\tilde r =r_*- R \,\, \mbox{ and} \,\,\, \tilde R=R-\e$$ 
  and suppose that  $ R<4\e<r/2$ so that 
  $${\textstyle  |r-r_*|<\e < \frac18 r, \,\,\, |r-\tilde r| < \frac14 r, \,\,\, \tilde R< 3\e \,\, \,\mbox{ and} \,\,\, \tilde r >\frac14 t.}$$ 
  We apply,    for  $s>\e(R+\e)$,  the first inequality of the  lemma that we have already proved at the beginning of this proof and,  for  $s \leq\e(R+\e)$,  Lemma \ref{lem3.4} with $\la=3$  
to infer  that 
  $$ \sup_{\xi\in \partial U(R)} \frac{h^*_\e(\xi, s)}{\e^{2\nu+1}} \leq \k_d\bigg(\frac1{\e} \vee \frac{\tilde R}{s}\bigg)p_s^{(d)}(\tilde R).$$
Thus  the repeated integral in (\ref{951})  is dominated by  a constant multiple of 
  \[
I := \int_0^{t}  \bigg(\frac1{\e} \vee \frac{\tilde R}{s}\bigg) p_s^{(d)}(\tilde R) q_{R}(r_*, t-s)ds. 
   \]  
   
  Write $I_{[a,b]}$ for the integral above  restricted on the interval $[a,b]$.    Applying  Theorem B  we see 
\[
I_{[0,t/2]}\leq   \k'_d\int_0^{t/2}   \bigg(\frac1{\e} \vee \frac{\tilde R}{s}\bigg)   
p_s^{(d)}(\tilde R)  \frac{\tilde r}{t-s}p_{t-s}^{(1)}(\tilde r) \bigg(\frac{R}{r}\bigg)^{(d-1)/2} ds\Big[1+ O\Big(\frac{t}{Rr}\Big)\Big].
  \]
On using the inequality $1/(t-s) \geq 1/t + s/t^2$,  the right-hand side  is bounded above by 
\beq
  \frac{\k''_d\,\tilde r}{t^{3/2}} \bigg(\frac{R}{r}\bigg)^{(d-1)/2}  e^{-\tilde r^2/2t}\int_0^{\infty} \bigg(\frac1{\e} \vee \frac{\tilde R}{s}\bigg)\exp\Big\{-\frac{\tilde r^2s}{2t^2}- \frac{\tilde R^2}{2s}\Big\}\frac{ds}{s^{d/2}} \Big[1\vee \frac{t}{\e r}\Big].
\eeq
Supposing
\beqn\label{000}
\tilde R \tilde r/t >1/2,
\eeqn
we compute the last integral  
(use if necessary  (\ref{13}) 
of Section {\bf5.2})  to conclude
$$I_{[0,t/2]}\leq \k_d'''\bigg(\frac1{\e} \vee \frac{\tilde r}{t}\bigg)\frac{1}{t^{d/2}}\bigg(\frac{R}{\tilde R}\bigg)^{(d-1)/2}  e^{- (\tilde r+ \tilde R)^2/2t}e^{\tilde R^2/2t}\Big[1\vee \frac{t}{\e r}\Big].
$$
For the other interval $[t/2,t]$ we obtain 
$$\bigg(\frac1{\e} \vee \frac{\tilde R}{t}\bigg)^{-1}I_{[t/2,t]} \leq  \frac{\k_d}{t^{d/2}}\int_{0}^{t/2} q_R(r_*,s)ds \leq  \frac{\k_d}{  t^{d/2}}P^{BM}_0\Big[\max_{s\leq t/2}X_s > r_*-R\Big],$$
and, since the last probability is at most $2e^{-2(r_*-R)^2/t}$,
taking  $R= 2 \e$  (so that $\tilde R= \e$ and $\tilde r +\tilde R =r_*-\e$) yields
$$I_{[t/2,t]} \leq \k_d'  \bigg(\frac1{\e} \vee  \frac{\e}{t}\bigg)\,p^{(d)}_{t/2}(r_*-2\e),$$
which combined with the bound of $I_{[0,t/2]} $ obtained above shows 
 $$I \leq \k_d'''   \bigg(\frac1{\e} \vee \frac{r}{t}\bigg)\, p_t^{(d)}(r_*-\e) e^{\e^2/2t}\Big[1\vee \frac{t}{\e r}\Big].
$$
provided  $r/t>1$ and (\ref{000}) is true. We may suppose $r^2> 8t$. For  if $r^2\leq  8t$,  entailing $r<8$ and $p_t^{(d)}(r) \asymp p_t^{(d)}(0)$, the formula to be shown follows from Lemma \ref{lem3.4} with $\la =8$.  Now take  $\e=t/r$, which conforms to the requirement (\ref{000})  as well as  the condition  $\e< r/8$ imposed at the beginning of the proof.  Then, 
   $p_t^{(d)}(r_*-\e) e^{\e^2/2t}\leq p_t^{(d)}(r-2\e) e^{\e^2/2t} \leq  p_t^{(d)}(r)e^{2\e r/t} = p_t^{(d)}(r)e^{2}$, and we  find that   $h^*_1(\z,t) \leq \k_d (r/t)p_t^{(d)}(r)$ as asserted in the lemma.
\qed 

\v2\n
{\sc Remark 5.}  The inequality in (\ref{951}) 
though appearing   intuitively  obvious may require verification.  We suppose  $d=2$ for simplicity 
and use the notation  $h_a^*(\z,t,\th)$, $-\pi < \th <\pi$, given  in (\ref{h*}) of the next section (it designates the density of $(\sigma_a, {\rm Arg} \, B(\sigma_a))$ at $(t,\th)$).  Write    $0'$ for  $(1-\e)\be$. 
For any $1< b< z$, the Brownian motion starting at $\z$  hits $\partial U(b)$ before $U(\e)+0'$ (the shift of $U(\e)$ by $0'$),   hence   
\beqn\label{Rem5}
h_\e^*(\z-(1-\e)\be,t)  = \frac{1}{2\pi}\int_{-\pi}^{\pi} d\phi \int_{0}^t h^*_{b}(\z,t-s, \phi) 
h_\e^*(be^{i\phi}-0',s)ds.
\eeqn
By using  an explicit form of the Poisson kernel of the domain  $\C\setminus U(\e)$ we deduce that for each  $\de>0$ (chosen small),  as $y:=b-1 \downarrow 0$ and  $\phi \to 0$
\beqn\label{Rem51}
 \frac{1}{2\pi \e}\int_0^\de h_\e^*(be^{i\phi}-0',s)ds = \frac1{\pi}\cdot \frac{y}{y^2+ (b\phi)^2}(1+o(1))
 \eeqn
 (cf. Appendix (C)). Restricting the range of the outer integral  to $|\phi|< \sqrt y$ in (\ref{Rem5})  and passing to the limit  we obtain the required upper  bound of $h_1^*(\z,t,0) =h_1^*(\z,t)$.


  \v2\n 
  {\bf 4.2.}  {\sc   Refinement in Case  $t<1$.} \,  
  In the next section we shall apply Lemma \ref{lem3.4} with  $\z$ on the  plane that is tangent to $U(a)$ at a point of the surface $ \partial U(a)$.  By the underlying  rotational invariance   we may suppose  that the plane is tangent at  $a\be$ so that  $\z\cdot \be = a$. Let $\phi$ be  the colatitude of $\z$ so that 
 \beqn\label{eta_y}\eta:=|\z- a\be|= a\tan \phi  \quad\mbox{and}\quad y:= |\z|-a= a\sec \phi -a.
 \eeqn
    Then  $y/a\sim \frac12 \phi^2$ and an elementary computation yields 
   \beqn\label{Texp}
   \phi^2+ \frac{y^2}{a^2} = \phi^2 + (\sec \phi -1) ^2 = \frac{\eta^2}{a^2} - \frac{5}{12}\phi^4 - O(\phi^6),
   \eeqn
  from which one may infer in one way or another that 
  in the case when  $y/\sqrt t$ is large the upper estimate of Lemma \ref{lem3.4} is  not fine enough: in fact   the term  $-\frac5{12}a^2\phi^4/2t$ can  be removed from the exponent of the exponential factor involved in  $p_t^{(1)}(y) p_t^{(d-1)}(a\phi)$ as asserted in the next proposition (cf. its Corollary).  (However, the bound of Lemma \ref{lem3.4} is  of correct  order  if  $\sqrt t > y$.)  This seemingly minor flaw  becomes serious in the proof of Theorem \ref{thm1.2} (when $\th$ is close to $\frac12 \pi$). 
  
 The next proposition partially  improves both Lemma \ref{lem3.4}  and the second inequality of Lemma 
 \ref{lem3.5}  (in the case   $t<1$).  Remember   the definition of   $h^*_a(\z,t)$ given  right after  (\ref{h_z0}). 
  \v2
\begin{Prop}\label{cor_Imp}
Let $z=|\z|>a$,   $\z\cdot \be/z = \cos \phi$ $(|\phi|< \pi)$, $y= z-a$,  and $r= |\z- a\be|$ as in Lemmas \ref{lem3.4} 
and \ref{lem3.5}.  There exist   positive constants  $C_1$, $C_2$ and   $C$ depending only on $d$ such that  whenever     $ t<a^2,   y<a$ and  $ |\phi| <1$,
 \beq
 \frac{C_1 a^{2\nu+1}y}{t}p_t^{(d)}(r) e^{-C[a\phi)^2(\phi^2+a^{-1}\sqrt t]/ t} \,\leq  \, h^*_a(\z,t)  \,
 \leq  \, \frac{C_2 a^{2\nu+1}y}{t}p_t^{(d)}(r)
 e^{C\phi^4[ay+(a\phi)^2]/t}.
 \eeq
 \end{Prop}
\v2

\begin{Cor} \label{prop_main}  There exists positive constants  $\k_d$ and  $M,$  depending only on $d$ such that if $t<a^2$ and $\eta$, $\phi$ and $y$ are given as  in (\ref{eta_y}) with  $|\phi| <1$, then 
$$
  h_a(a+y,t, \phi) \leq \frac{\k_d a^{2\nu+1} y}{t}p^{(d)}_t(\eta) e^{M\eta^6/a^4t }.
  $$
    \end{Cor}
 \v2
 
Our  proof  of Proposition \ref{cor_Imp} rests on the skew product formula (\ref{h_skew}) and requires  some elaborate  estimate of the distribution  of 
the  random time $\tau^X$  given in (\ref{tau_1}), namely
$$\tau^X = \int_0^{T_0} \frac{ds}{ (1+ X_s)^2}.$$
Here  $X$ denotes a standard linear Brownian motion; its law conditioned on $X_0=r$ is denoted by  $P_r^X$ as mentioned in the beginning of this section. 
 In the situation we are interested in,   the starting point of $B_t$ is close to the sphere $\partial U(1)$, so that the Bessel process $|B_t|$  may be replaced by linear Brownian motion. 

\begin{Lem}\label{claim} \, For $b>0$ and $r>0$,
\beqn\label{trunc} 
P^{BM}_r[ X_{1-s} \geq br +rs \,\,\mbox{for some}\,\, s\in [0,1]\,|\, T_0=1] \leq 6 e^{- \frac23 b^2r^2}.
\eeqn
\end{Lem}
\v2\n
\pf\,
 Let  $R_t, t\geq 0$ be a three-dimensional Bessel process and $L_r$ its last passage time of $r$. Then we have the following sequence of identities of conditional laws:
\begin{eqnarray}\label{id_law} 
&&(X_{1-s})_{0\leq s\leq 1} \,\,\mbox{conditioned on}\,\,X_0=r, T_0=1\nonumber\\
 &&\quad \stackrel{\rm law}{=} (R_s)_{0\leq s\leq 1} \,\,\mbox{conditioned on}\,\, R_0=0, L_r=1 \nonumber\\
 &&\quad \stackrel{\rm law}{=} (R_s)_{0\leq s\leq 1} \,\,\mbox{conditioned on}\,\, R_0=0, R_1=r\\
 &&\quad \stackrel{\rm law}{=} (R_{1-s})_{0\leq s\leq 1} \,\,\mbox{conditioned on}\,\, R_0=r, R_1=0   \nonumber\\
&&\quad \stackrel{\rm law}{=} (sR_{s^{-1}-1})_{0\leq s\leq 1} \,\,\mbox{conditioned on}\,\, R_0=r \nonumber
 \end{eqnarray}
 (see \S 1.6 and \S 8.1   of \cite{YY} and (3.7) and (3.6) in \S XI.3 of \cite{RY}). On using the last expression  a simple manipulation shows that the conditional probability in (\ref{trunc}) equals 
\beqn\label{EQ_r}
P^R[ R_u > (b+1)r+ bru   \,\,\mbox{for some}\,\, u \geq 0\,|\, R_0=r ],
\eeqn
where $P^R$ denotes the law of $(R_t)$.
 Since $R_t$ has the same law as the distance from the origin  of a three-dimensional Brownian motion starting at $(r/\sqrt 3, r/\sqrt 3,r/\sqrt 3)$, the  probability in (\ref{EQ_r})   is dominated by 
 $$3P^{BM}_{r/\sqrt 3}\Big[ |X_s| > (b+1+ bs)r/\sqrt3   \,\,\mbox{for some}\,\, s \geq 0 \Big],$$
which is at most $6 e^{- \frac23 b^2r^2}$ according to a well known bound of escape probability  of a linear Brownian motion with drift. The bound (\ref{trunc}) has been verified.  \qed

\begin{Lem}\label{Imp} \, There exists a  constant $C>1$ such that for $0< \de \leq 1$, $0<t<1$ and $y > 0$, 

 ${\rm (i)} \quad {\displaystyle P^{BM}_y\Big[\tau^X \geq \frac{t}{1+ (1-\de)y}\,\Big|\, T_0 =t\Big]\leq  
C\bigg(1\wedge \frac{\sqrt t}{\de y}\bigg)e^{-3[\de (1- 2y)]^2y^2/2t}\quad \mbox{if}\quad y<\frac14.}$
\vskip3mm 
 ${\rm (ii)} \quad {\displaystyle P^{BM}_y\Big[\tau^X \geq \frac{t}{1+ (1+\de)y+ \de y^2}\,\Big|\, T_0 =t\Big]\geq  1- C^{-1} e^{-\de^2y^2/6t}.}$
\end{Lem}
  \v2\n
  \pf\, By the scaling property  of  $X$ the conditional probabilities to be estimated may be written as
  \v2
  $I_{-}:= P^{BM}_r[\tilde\tau^X \geq \frac{1}{1+ (1-\de)y}\,|\, T_0 =1]\quad$ and 
  $\quad I_{+}:= P^{BM}_r [\tilde\tau^X \geq \frac{1}{1+ (1+\de)y+ \de y^2}\,|\, T_0 =1],$
  \v2\n
  where
  $$ r=\frac{y}{\sqrt t},\,\, \tilde \tau^X = \int_0^1 \frac{ds}{(1+ \sqrt t X_s)^2}.$$
  
According to Lemma \ref{claim}   the lower bound (ii) readily follows  from  this expression. Indeed,   if $\sqrt t X_{1-s} < ys +\frac12 \de y$ for $0<s<1$, then
$$ \tilde \tau^X \geq \int_0^1\frac {ds}{(1+ys+\frac12\de y)^2}= \frac1{1+(1+\de )y+(1+\frac12\de)\frac12 \de y^2},$$
implying the occurrence of the event of the conditional probability giving  $I_+$, hence the required lower bound.

 The upper bound (i)  requires a  delicate estimation. We write the event under the conditional probability  for $I_-$ in the form 
\begin{eqnarray}\label{write}
\tilde\tau^X - \frac{1}{1+ y} &=& \int_0^1\bigg[\frac1{(1+\sqrt t X_s)^2} -\frac1{(1+ys)^2}\bigg]ds  \\
&\geq& \frac{\de y}{(1+(1-\de)y)(1+y)}.\nonumber
\end{eqnarray} 
 Observe that the integral above is less than $2\int_0^1(ys-\sqrt t\, X_s)ds$ a.s.  and the last member is larger than $\de y(1-2 y)$ (for $y>0$), so that the inequality   (\ref{write}) implies 
  \beqn\label{write2}
\int_0^1(ys-\sqrt t\, X_s)ds \geq 
\frac12 \de y(1  -  2y) \quad \mbox{if} \quad \sup_{0<s<1}|X_s -rs |< 2r.
\eeqn
Owing  to Lemma \ref{claim}  we have $P^{BM}_0[ \sup_{0<s<1}|X_s -rs |>2r \,|\, T_r=1] \leq 12 e^{- 2y^2/t} $, 
  which along with (\ref{write2})  shows
  $$I_-\leq P^{BM}_0\bigg[\int_0^1(ys-\sqrt t X_s)ds \geq  \frac12 \de y(1  -  2y)\,\bigg|\, T_r =1\bigg] + 12e^{-2y^2/t}.$$
 Using (\ref{id_law}) again we rewrite the probability on the right  in terms of the three dimensional Bessel process $R_t$, which results in 
$$P^R\bigg[\int_0^\infty \frac{r-R_s}{(1+s)^3}ds>\frac12\de r(1 - 2 y) \,\bigg|\, R_0 =r\bigg].  $$ 
For our present objective of obtaining an upper bound we may replace $R_s$ by $X_s$. Since the random variable 
$\int_0^\infty \frac{r-X_s}{(1+s)^3}ds= \frac12\int_0^\infty (1+s)^{-2}dX_s$   is Gaussian of mean zero under $P^{BM}_r$ and  its variance equals
$$ E^{BM}_0\bigg[\Big(\int_0^\infty \frac{ X_s ds}{(1+s)^{3}}\Big)^2\bigg]=   \frac14\int_0^\infty (1+s)^{-4}sds  = \frac1{12}, $$ 
it follows that if  $y<1/4$,
$$I_- \leq  C\bigg(1\wedge \frac{1}{\de r}\bigg)e^{-3r^2(\de - 2\de y)^2/2} +12e^{-2y^2/t}.$$
On the right-hand  side the second term may be absorbed into the first,  resulting in  the required bound.
\qed

\v2\v2
The next lemma, valid for all $d\geq 2$,  improves  the bound of Lemma  \ref{lem3.4} when  $r/t>1$.  
  \v2\n
  \begin{Lem}\label{Imp1} \,   There exists a positive constant  $C$  depending only on $d$ such that 
  \beqn\label{crucial}
  h_a(a+y,t, \phi) \leq \frac{Ca^{2\nu+1}y}{t^{1+d/2}} 
  \exp\Big\{- \frac{1}{2t}\Big( (a^2+ ay)\phi^2 + y^2 -\frac{a^2}{12} \phi^4 - 12a y\phi^4 \Big) \Big\}. 
  \eeqn
 whenever     $0<  y <a$, $t<a^2$ and $|\phi|<1$.
\end{Lem}
  \v2\n
  \pf\,  Suppose   $d=2$ and  $a=1$, the case   $d\geq 3$ being briefly  discussed at the end of this proof.  Let $\tau$  be as in the preceding lemma. From  (\ref{h_skew2})  it plainly follows that
    \beqn\label{h_skew3}
h_1(1+y,t, \phi)\leq 2\pi E^{BM}_y[e^{\tau^X/8}p^{(1)}_{\tau^X}(\phi) \,|\, T_0=t]q^{(1)}_1(1+y,t).
   \eeqn
 Noting   $\tau^X< t$,  we compute   $E_y^{BM} [e^{-\phi^2/2\tau^X}\,|\, T_0= t]$.    
    Define the random variable $\mit\Delta$ via
$$\frac1{ \tau^X} = \frac{1+y - y\mit\Delta}{t},$$ so that
\beqn \label{so that}
E_y^{BM}[e^{-\phi^2/2\tau^X}| T_0= t] = e^{-(1+y)\phi^2/2t}E_y^{BM}[e^{(\phi^2/2t)y\mit\Delta}| T_0= t]. 
\eeqn
Put  
$F(\de)=E_y^{BM}[\mit\Delta \geq \de$ $|\, T_0= t ]$ for $  -\infty<\de\leq  1$. 
Then by Lemma \ref{Imp} (i)
$$F(\de) = P^{BM}_y\Big[ \tau^X \geq \frac{t}{1+ (1-\de)y}\,\Big|\, T_0 =t\Big]\leq \frac{ C}{1+\de y t^{-1/2}\,} e^{-3[\de(1 - 2y)]^2 y^2/2t}$$
(for $y<1/4, 0<\de\leq 1$).    Put 
$$A = \frac{\phi^2}{2t}y \quad\, \mbox{and} \quad\, B= A\frac{\phi^2}{y} = \frac{\phi^4}{2t}.$$
 Then,    noting $F(1-0)=0$, we  perform  integration  by parts to see that
\beq
E_y^{BM}[e^{(\phi^2/2t)y\mit\Delta}| T_0= t] &=& -\int_{-\infty}^1 e^{A\de}dF(\de) = \int_{-\infty}^1 A e^{A\de}F(\de) d\de\\
&\leq& 1+ C\int_{0}^1 \frac{ A}{1+\de y t^{-1/2}\,} \exp\Big\{A\de - 3 \frac{\de^2(1-2y)^2y^2}{2t}\Big\}d\de.
\eeq

The last integral  restricted to the
 interval  $(\phi^2/y)\wedge 1\leq \de \leq 1$ is dominated  by 4,  provided  $y<1/8$, for  in this interval we have $\de y \geq \phi^2$ so that  the exponent  involved  in  the integrand is bounded from above by
 $$A\de - 3\frac{\de^2y^2}{2t}  \Big(1 - \frac14\Big)^2  \leq  - \frac14 A\de$$
  ($y\leq 1/8$), and thus the integral by  $\int_0^\infty A e^{-A\de/4}d\de =4$. 
On the other hand,   write the exponent  as
$$A\de - 3 (1-2y)^2\frac{\de^2y^2}{2t} = \frac{B}{12} - 3\bigg(\frac{ y}{\sqrt{2t}}\de - \frac{1}{6}\sqrt B\bigg)^2 + 3\frac{4 \de^2 y^3(1 -y)}{2t}$$
and  observe  that  the last term is less than  $6\phi^4 y/t$ if $\de\leq \phi^2/y$. Then,  we transform    the integral over $[0, \phi^2/y)$ by  changing the variable of integration according to   $u=\frac{ y}{\sqrt{2t}}\de - \frac16 \sqrt B$ and,  noting   $A\sqrt{2t}/y = \sqrt B$,   we deduce that  it  is  at most
 $$\sqrt B \int_{-\sqrt{B}/6}^{5\sqrt{B}/6}\frac{e^{-3 u^2}du}{ 1+ \sqrt 2 (u + \frac16\sqrt{B})}  \exp\Big\{\frac{B}{12} + \frac{6y\phi^4}{t}\Big\}\leq \frac{C\sqrt{B}}{1+\sqrt B}\exp\Big\{\frac{B}{12}  + \frac{6y\phi^4}{t}\Big\},
 $$
hence by virtue  of (\ref{so that})
\begin{eqnarray}\label{exp9} 
E_y^{BM}[e^{-\phi^2/2\tau^X}| T_0= t] &\leq& Ce^{-(1+y)\phi^2/2t}\bigg(1 + \frac{\sqrt{B}}{1+\sqrt B} \exp\Big\{ \frac{ \frac{1}{12} \phi^4 + 12 y\phi^4 }{2t} \Big\}\bigg) \nonumber \\
&\leq& C'  \exp\Big\{ \frac{ -(1+y)\phi^2 + \frac{1}{12} \phi^4 + 12 y\phi^4 }{2t} \Big\}.
\end{eqnarray}
On recalling  (\ref{h_skew3}) (and  (\ref{1dim}) as to $q_1^{(1)}$)   this concludes  
 the assertion of the  lemma, for  if $\phi^2 > t$, then $p^{(1)}_{\tau^X}(\phi) \leq p^{(1)}_{t/2}(\phi)$ on the event $\tau^X\leq t$ and $p^{(1)}_{\tau^X}(\phi)\leq \frac1{\sqrt {2\pi t}}e^{-\phi^2/2\tau^X}$ oherwise, while if $\phi^2\leq t/2$,  the lemma is obvious (see e.g.  Lemma \ref{lem3.4}).

The  case $d\geq 3$ is dealt with in the same way as above for the same reason mentioned at the end of the proof of Lemma \ref{lem3.4}.  We employ  the skew product representation of $d$-dimensional Brownian motion. For the radial component the same remark as given in Remark 4 is applied to the bounds obtained in Lemma \ref{Imp}. The spherical component behaves as the Brownian motion on the flat space for small $t$. It follows that in place of (\ref{h_skew2})  we have
$$h_1(1+y,t, \phi)  \leq  C_d E^{BM}_y[p^{(d-1)}_{\tau^X}(\phi)\,|\, T_0=t]q^{(1)}_1(1+y,t).$$
Thus 
the desired bound (\ref{crucial}) follows from  (\ref{exp9}).
 \qed


\v2\n
{\it Proof of Proposition \ref{cor_Imp}.}\,  For   $|\phi|<1$,
\begin{eqnarray}
 |\z- a\be|^2 &=& (y+a)^2 - 2(ay+a^2)\cos \phi +a^2 \nonumber\\
&=&y^2+  (a+y)\phi^2 -\frac{a^2}{12}\phi^4 + O(\phi^4 ay + a^2\phi^6),
\label{cos_f}
\end{eqnarray}
and the upper bound in Proposition \ref{cor_Imp} follows from  Lemma \ref{Imp1}.
  
 For the lower bound   we suppose  $a=1$ for simplicity and  apply  the skew product expression (\ref{h_skew}).  Suppose $d=2$. As in the proof of  Lemma \ref{lem3.4} (see (\ref{h_skew2})) we have
 $$
h^*_1(\z,t) = h_1(1+y,t,\phi) \geq E_y^{BM}[p_{\tau^X}^{(1)}(\phi) \,|\, T_0=t] q^{(1)}_0(y,t).
 $$ 
Plainly    $\tau^X<t$ from the  definition, while by   (ii) of Lemma \ref{Imp} with $\de = \sqrt t \,/y$ we  see that
 $P^{BM}_y[\tau^X>t (1+y + \sqrt t(1+y))^{-1}\,|\, T_0=t] \geq 1-e^{-1/6}$.   If  $t < \phi^2$ (so that  $p_\tau(\phi)$ is increasing in $\tau \in (0, t]$), it therefore follows that
 $$h^*_1(\z,t) \geq \k_d \frac{y}{t^2} \exp\bigg\{- \frac{y^2 + (1+y)\phi^2+ 2\phi^2\sqrt t\, )}{2t}\bigg\},$$
 entailing the required lower bound in view of (\ref{cos_f}).
If $t\geq \phi^2$, then    the conditional expectation above  is bounded  below by 
 $\k_d/\sqrt t$ and observing  $(r^2-y^2) /t \leq (1+y)$ we obtain $h^*_1(\z,t) \geq \k_d'yt^{-1}p^{(2)}_t(r)$, a better lower bound.
  \qed

\section{ Proof of Theorem \ref{thm1.2}  (Case  $d=2$) }
  Throughout this section we let   $d=2$; also  let  $\x =x{\bf e}$ and write  $v$ for $x/t$. The definition of $h_a$ given at the beginning  of Section 4 may read 
$$h_a(x,t, |\th|) = 2\pi P_\x[{\rm Arg}\,  B(\sigma_a)\in d\th, \sigma_a\in dt]/d\th dt \quad (x>a, -\pi <\th <\pi).$$
In this section  we prove 
\vskip20mm
\begin{Thm}\label{thm5.1}  Let $v=x/t$. Then, 
\v2
{\rm (i)} uniformly  for $ 0\leq \th < \frac12\pi - v^{-1/3} $  and  for  $t>1$, as $v\to \infty$
\beqn\label{EQ}
h_a( x,t, \th) =  2\pi av \, p_t^{(2)}(|\x -ae^{i\th}|) \cos \th \bigg[
1 + O\Big( \frac{1}{(\frac12\pi -\th)^3 v}\Big) \bigg];\, \mbox{and}
\eeqn

{\rm  (ii)} there exists a universal   constant $C$  such that  for $ |\frac12\pi -\th|<(av)^{-1/3}$, $t>a^2$ and $v>2/a$, 
$$C^{-1} \frac{1}{(av)^{1/3}}\leq \frac{h_a(x, t, \th)}{a v e^{-av(1-\cos \th)} p_t^{(2)}(x -a)} \leq C \frac{1}{(av)^{1/3}}. $$
\end{Thm}

\v2
Note that $\cos \th \sim \frac12\pi -\th$ as $\th \to \frac12\pi$ and 
\beqn\label{eq5.10}
p_t^{(2)}(|\x -ae^{i\th}|) =e^{-av(1-\cos \th)} p_t^{(2)}(x -a).
\eeqn
The following corollary of Theorem \ref{thm5.1}  is  a restatement of Theorem \ref{thm1.2} for the case $d=2$.

\begin{Cor}\label{cor1}   For  $t>a^2$ and  $v= x/t>2/a$, 
\beq
 &&\frac{P_\x[{\rm Arg}\, B(\sigma_a)\in d\th\,|\, \sigma_a =t]}{d\th}  \\
 &&\quad = \sqrt{\frac{av}{2\pi}}e^{-av(1-\cos \th)}\cos \th \bigg[ 1 + O\Big(\frac{1}{av\cos^3 \th}\Big)\bigg]   \quad \mbox{if}  \quad \cos \th \geq\frac{1}{(av)^{1/3}};\, and  \\
 &&\quad \asymp   \sqrt{av}\,e^{-av(1-\cos \th)}  (av)^{-1/3}  \quadd\quadd\,\,   \mbox{if} \quad   |\cos \th| \leq \frac{1}{(av)^{1/3}}.
\eeq 
\end{Cor}

For $|\th| > \frac12 \pi + (av)^{1/3}$ we shall obtain  an upper bound  (Lemma \ref{wc}) which together with Corollary \ref{cor1}  verifies the next corollary.
\begin{Cor}\label{cor2} 
As $v:=x/t\to\infty$ 
$$  \sqrt{\frac{\pi}{2av}}e^{av(1-\cos \th)} P_{\x}[{\rm Arg}\, B(\sigma_a)\in d\th\,|\, \sigma_a =t] \,\Longrightarrow \,\frac{1}{2}\1\Big( |\th| <\frac12 \pi\Big) \cos \th \,d\th.$$
\end{Cor}

\v2
For the proof it will become convenient to bring in the notation
\beqn\label{h*}
h^*_a(\z,t, \th) =2\pi \frac{P_\z[{\rm Arg}\,  B(\sigma_a)\in d\th, \sigma_a\in dt]}{d\th dt} \quad (|\z| >a, 0 \leq |\th |< \pi).
\eeqn
which  is a natural  extension of $h^*_a$ introduced in Section {\bf 4.1}:  $h^*_a(\z,t)= h^*_a(\z,t,0)$; also,  
$h_a(z, t, |\th|) = h^*_a(z\be,t, \th)$.

 \v2\n
{\bf 5.1.} {\sc Lower Bound I}. 
\v2

The following lemma, though easy to obtain, gives a correct asymptotic form of $h_a$ if $\th \in (0, \pi/2)$ is away from $\frac12 \pi$  and provides a guideline for later arguments.  Combined with Theorem A  it also entails
Proposition \ref{thm0}. Let $\x =x\be$ and $v=x/t$ and put
$$\Psi_a(x,t,\th)  = \frac{2\pi ax}{t} e^{-\frac{ax}{t}(1-\cos \th)} p_t^{(2)}(x -a) \Big(\cos \th -\frac{a}{x}\Big).
 $$

\v2
\begin{Lem}\label{LBD}  \, For all $x>a,  t>a^2$ and $\th\in (-\frac12 \pi,\frac12 \pi)$, 
\beqn\label{LB}
 \frac{P_\x[\arg B(\sigma_a)\in d\th, \sigma_a\in dt] }{d\th dt} \geq  \Psi_a(x,t,\th);
 \eeqn
 in particular 
$h_a(x,t, \th) \geq  \Psi_a(x,t,\th)$.
\end{Lem}

 \begin{figure}[b]
 \begin{center}
\begin{picture}(460,145)
\put(36,110){\line (1,0){374}}
\put(408,140){\line (0,-1){126}}

\put(408,38){\line(-4,1){394}}
\put(390,43){\vector (-4,1){40}}
\put(391,42){\vector (1, 4){10}}
\put(391,42){\line (1, 4){24}}
\put(41,130){\vector(4,-1){30}}
\put(408,17){\line (-4,1){373}}
\put(36,110){\line (1,4){5}}

\put(36,110){\circle{40}}
\put(36,110){\circle*{4}}
\put(41,130){\circle*{3}}
\put(408,110){\circle*{4}}
\put(391,42){\circle*{3}}

\put(26,106){$0$}
\put(40,113){$\th$}
\put(398,22){$\th$}
\put(398,114){$\x$}
\put(62,129){$\eta$}
\put(40,135){$ae^{i\th}$}
\put(390,79){$u$}
\put(364,53){$l$}
\put(240,86){$L(\th)$}
\put(184,8){Figure 1}
\end{picture}
\end{center}
\vspace*{0cm}
\end{figure}

\v2\n
\pf\, We represent  points on the plane by complex numbers. Let $0\leq \th <\frac12 \pi$ and denote by $L(\th)$ the straight   line tangent  to the circle $\partial U(a)$ at $ae^{i\th}$.   Let $\sigma_{L(\th)}$ be the first time $B_t$ hits  $L(\th)$ and consider the coordinate system  $(u, l)$ where  the  $u$-axis is the line through $\x$ perpendicular to $L(\th)$  and  the $l$-axis is $L(\th)$ so that the $l$-coordinate of the tangential point $ae^{i\th}$ equals $x\sin \th$ (see Figure 1).  Put 
\beqn\label{psi0} 
\psi_a(l,t)= \frac{P_\x[ B(\sigma_{L(\th)})\in d l, \sigma_{L(\th)} \in dt]\,}{dl dt}
\eeqn
and 
 \beqn\label{K}
U = \int_0^{t} ds \int_{\R\setminus \{x\sin \th\}} \psi_a(l, t-s)
h^*_a(\xi^*_a(l),s,\th)dl,
\eeqn
where  $h^*_a$ is defined by (\ref{h*}) and  $\xi_a^*(l)$ denotes the point of the plane which lies  on $L(\th)$ and whose $l$-coordinate equals $l$ (so that $\xi^*_a(x\sin \th) = ae^{i\th}$).
  Then
  \beqn\label{LB9}
h_a(x, t, \th) = h^*_a(\x,t,\th) = 2\pi a \psi_a(x\sin \th,t)   + U.
 \eeqn
 Here  the factor $a$ of the first term  on the right-hand side of (\ref{LB9})  comes out from  the relation $dl = ad\th$ valid at $ae^{i\th}$; 
  for the present proof  we need only the lower bound (with $U$ being discarded)  that is verified  by the same argument as in Remark 5; the equality however is used later,  whose verification we  give after  this proof.  We claim 
\beqn\label{clm8}
  \psi_a(x\sin \th,t)  =  \frac1{2\pi a}\Psi_a(x,t,\th).
 \eeqn
Since the $u$-coordinate of $\x$ equals  $x\cos \th -a$   we have in turn
  $$
 \psi_a(l,t) =  \frac{x\cos \th -a}{t} \, p^{(1)}_{t}(x\cos \th -a)p^{(1)}_{t}(l)$$ 
and
\beqn\label{eq3.1}
  \psi_a(x\sin \th,t) 
= \frac{x\cos \th -a}{2\pi t^2} e^{- |x e^{i\th} -a|^2/2t }.
\eeqn
Hence, noting (\ref{eq5.10}),   we  readily identify the right-hand side of (\ref{eq3.1})   with that of (\ref{clm8}).

Finally one may realize that  (\ref{LB9})  shows    $a\psi_a(x\sin \th,t)$  to be  a  lower bound     for  the density  of the distribution of $(\arg B(\sigma_a), \sigma_a)$ (rather than  $({\rm Arg}\, B(\sigma_a), \sigma_a)$).   \qed
\v2

 {\it Proof of (\ref{LB9}).}  
We are to take the limit as  $b\downarrow a$ in  the expression
 \beqn\label{m9}
h_a^*(\x,t,\th) = \int_0^tds \int_{l\in \R} \psi_b(l,t-s)h^*_a(\xi^*_b(l), s,\th)dl \qquad (a<b <x).
 \eeqn
 Here the coordinate $l$ and  $\xi^*_b(l)$ are analogously defined with the tangential line to $\partial U(b)$  at $be^{i\th}$.
 Put $y=b-a$ and split  the inner  integral in (\ref{m9}) at $l=   x\sin \th \pm \sqrt y$.

First consider 
 $$
m_{\rm in}(b) := \int_0^tds \int_{|l- x\sin \th|<\sqrt y} \psi_b(l,t-s)h^*_a(\xi^*_b(l), s,\th)dl.
$$
 As in Remark 5 we apply  the explicit form of the Poisson kernel of   $\C \setminus  U(a)$ to see  that for each $\de>0$,  uniformly for $l: |l-x\sin \th|<\sqrt y$, as $y\downarrow 0$
  $$ ({2\pi a})^{-1}\int_0^\de h^*_a(\xi^*_b(l), s,\th) ds =\frac{1}{\pi}\cdot  \frac{y}{y^2+(l- x\sin \th)^2}(1+o(1),$$
 which yields  
 $
  \lim_{b\downarrow a} m_{\rm in}(b) = 2\pi  a\psi_a(x\sin\th, t)
$
in view of continuity of  $\psi_b(l,t-s)$.
    
  As for the contribution of the range $\{l: |l-x\sin \th|\geq \sqrt y \}$, denoted by  $m_{\rm out}(b)$, we 
  substitute  in the  integral representing  it  the expression    
  $$h^*_a(\xi^*_b(l), s,\th)=  \int_0^s ds'  f_{b-a}(l-l', s-s')\int_{l' \in \R}  h^*_a(\xi^*_a(l'), s',\th)dl'ds',
$$
where $f_{y}(l-l', s-s') = y (s-s')^{-1}p_{s-s'}^{(1)}(y)p_{s-s'}^{(1)}(l-l')$, representing the space-time hitting density of the line $L(\th)$ for the Brownian motion $B_t$ conditioned on  $B_s=\xi^*_b(l)$, and
 perform the integration w.r.t. $dsdl$ first to see that  $m_{\rm out}(b)$ converges to $U$.  \qed 
 \v2\n

{\bf 5.2.} {\sc Upper Bound I}.

\begin{Prop}\label{UBD1} \, Let $v=x/t>1$ and $t>1$.  For some  universal constant $C>0$,
$$h_a(x,t, \th) \leq  \Psi_a(x,t,\th)
\bigg[1+\frac{C}{ (\frac12 \pi -\th)^3 av} \bigg]\quad\quad \mbox{if}\quad 0 \leq \th <\frac{\pi}2 -\frac1{(av)^{1/3}} .
$$
\end{Prop}
\v2\n

For the proof of Proposition \ref{UBD1} we  compute     $U$  given in (\ref{K}): it suffices to  show  the  upper bound
\beqn\label{U}
U \leq \frac{C\Psi_a(x,t,\th)}{ (\frac12 \pi -\th)^3 av} \qquad \mbox{for}\quad 0 \leq \th <\frac{\pi}2 -\frac1{(av)^{1/3}}.
\eeqn
Let $\psi_a$ be the  density of  the hitting distribution  in space-time  of $L(\th)$ defined by (\ref{psi0}).   Bringing in the new variable $\eta\in \R$ by 
 $$l=  x\sin \th -\eta $$
 we write
 \beqn\label{psi01}
 \psi_a(l,t) =  \frac{x\cos \th -a}{t} \, p^{(1)}_{t}(x\cos \th -a)p^{(1)}_{t}(x\sin \th -\eta).
 \eeqn
  We break the repeated integral defining $U$ into  two parts by splitting the time interval  $[0,t]$ at $s=a/v$ (namely $s/a^2 = 1/av$, conforming to the scaling relation) and denote the corresponding integrals by  
 $$U_{[0,a/v] }\quad \mbox{ and}\quad  U_{[a/v,t]},$$
   respectively.
    The rest of the proof is divided   into three steps. 

 {\it Step 1.}  \, 
 The essential task  for the proof  is performed  in the estimation of  $U_{[0, 1/v]}$, which   is involved in  Lemmas  \ref{lem5.2.1} through \ref{lem5.2.4}.

 Recall  
  $$U_{[0,a/v]} = \int_0^{a/v} ds \int_{\R} \psi_a(l, t-s)h_a^*(\xi^*_a(l),s,\th)dl,$$
and  write
$$J_E = \frac1{a}\int_{E} e^{v\eta \sin \th } \,d\eta\int_0^{a/v} \exp\Big\{-\frac{v^2}{2} s\Big\}h_a^*(\xi^*_a(l),s, \th)ds \qquad (E \subset [0,\infty)).$$
 With an obvious 
 reason of comparison we may restrict our consideration to the half line  $l< x\sin \th$, i.e. to $\eta >0$.

 \begin{Lem}\label{lem5.2.1}  \quadd
 $U_{[0,a/v]} \leq   C\Psi_a(x,t,\th) J_{[0,\infty)}.$ 
 \end{Lem}
 \v2\n
 \pf\, 
We see from  (\ref{psi01})  
 $$
 \psi_a(l,t-s)  
 = \frac{x\cos \th -a}{t-s} e^{-ax(1-\cos \th)/(t-s)}\, p^{(2)}_{t-s}(x -a)\exp\Big\{\frac{2 x\eta \sin \th - \eta^2}{2(t-s)}\Big\}.
 $$ 
 On using $\frac1{t-s} = \frac1{t}+\frac{s}{t(t-s)}$  an elementary computation leads to
\begin{eqnarray}\label{psi00}
e^{- ax(1-\cos \th)/(t-s)}\, p^{(2)}_{t-s}(x -a) &=&  (1-s/t)^{-1}p^{(2)}_{t}(x-a) e^{- av(1-\cos \th)} e^{-v^2s/2}  \nonumber\\
&& \,\,\, \times \exp\Big\{\frac{-v^2s^2 +2avs \cos \th - a^2st^{-1}}{2(t-s)}\Big\}.
\end{eqnarray}
and substitution  in the preceding formula yields
\beq
\psi_a(l,t-s)
& =& \bigg(\frac{t}{t-s}\bigg)^2\frac1{2\pi a} \Psi_a(x,t,\th)e^{v\eta \sin \th}e^{-v^2s/2}\\
&& \,\,\times \exp\Big\{\frac{-(v^2s^2 +\eta^2 -2vs\eta \sin \th)  + 2avs \cos \th - a^2st^{-1}}{2(t-s)}\Big\}.
\eeq
With the help of the inequality  $v^2s^2 +\eta^2 -2vs\eta \sin \th  >0$ this   leads to
\beqn\label{Ineq3}
\psi_a(l,t-s) \leq \bigg(\frac{t}{t-s}\bigg)^2\frac1{2\pi a}\Psi_a(x,t,\th) e^{v|\eta| \sin \th }\exp\Big\{-\Big(\frac{v^2}{2} -\frac{av\cos \th}{t-s}\Big)s\Big\}
\eeqn
valid for all $0<s < t, |\eta|<\infty$.
 Now, in  (\ref{Ineq3}) we  get  a constant  to  dominate both the heading  factor and  the term $(av\cos \th)s/(t-s)$ in the exponent  for $s< a/v$  to  see  the inequality of the lemma. \qed
\v2

\v2
 {\it Step 2.}  \,  In this  step we prove three lemmas that together verify
  \beqn\label{U0}
 U_{[0,a/v]} \leq C \Psi_a(x,t,\th)\frac1{ av \cos^{3} \th} \quad\mbox{if}\quad \cos \th > \frac1{(av)^3} .
 \eeqn
  
 Let $\phi$ denote the angle between the rays $ra^{i\th}, r\geq 0$ and $r\xi^*_a(x\sin \th -\eta), r\geq 0$ so that
$$\eta = a\tan \phi \quad \mbox{ and} \quad y =a \sec \phi -a$$ 
and  $h_a^*(\xi^*_a(l),s, \th) = h_a(a+y,s, \phi)$.
Applying  Lemma \ref{lem3.4} with  $\la=\pi$, we infer  that for $0\leq b <b'\leq 1$, 
\beq
J_{[b,b']} &\leq & \frac{2\k_d}{a} \int_{b}^{b'} e^{v\eta \sin \th } \,d\eta \int_0^{a/v}   \frac{ay}{s} p_s^{(1)}(y) p_s^{(1)}(a\phi) e^{-v^2s/2}ds.
\eeq
Now and later we use the formula
\begin{eqnarray}\label{13}
\int_0^\infty \exp\Big\{-\frac{\eta^2}{2s}-\frac {v^2s}{2}\Big\}\frac{ds}{s^{p+1}} 
&=& 2\bigg(\frac{v}{\eta}\bigg)^p K_p(v\eta) \\
&\sim&  \left\{ \begin{array} {ll} 2^p\Ga(p) \eta^{-2p}\quad &  ( v\eta \downarrow 0, p>0) \nonumber\\[2mm]
{\displaystyle \bigg(\frac{v}{\eta}\bigg)^p\frac{ \sqrt{2\pi} \, e^{-v\eta} }{\sqrt{v\eta}} }\quad&  (v\eta \to \infty, p\geq 0)
\end{array} \right.
\label{K_p}
\end{eqnarray}
valid for all $\eta>0$ and $v>0$  (\cite{E}, p146). 
Noting  that since $y\sim \frac12 a\phi^2\sim \frac1{2a} \eta^2$,
\beqn\label{oo}
\frac{ay}{s} p_s^{(1)}(y) p_s^{(1)}(a\phi) \leq \frac{\eta^2}{s^{2}} e^{-(y^2+\phi^2)/2s} \qquad (\eta< 1),
\eeqn    
 we apply the equality in  (\ref{13}) with $\sqrt{y^2+\phi^2}$ in place of $\eta$  to deduce
\beqn\label{K0}
J_{[b,b']} 
\leq \frac{C}{a}\int_b^{b'}  \eta^2\frac{v}{\sqrt{y^2+\phi^2}} K_{1}(v\sqrt{y^2+\phi^2})  e^{v\eta  \sin\th} d\eta. 
\eeqn
Recall (\ref{Texp}), which may reduce to
\beqn\label{Ineq2}
y^2 + (a\phi)^2 > \eta^2(1 - Ca^{-2}\eta^2) \quad (|\phi|< 1)  
\eeqn
(for some $C>0$),
and   we evaluate 
  the  integral  over $\eta <a/v$ and  conclude the following 
    \begin{Lem}\label{lem5.2.2}
$$J_{[0,a/v]} \leq C\int_0^{a/v} e^{v\eta\sin \th}d\eta \asymp  \frac1v. $$
\end{Lem}

 In the rest of this proof of Proposition \ref{UBD1} we suppose   for simplicity 
 $$a=1.$$
 


\v2
The integral   $J_{[1/v, \infty)} $   may be easily evaluated with the same bound as above if $\th$ is supposed to be away from $\frac12 \pi$.  In order to include  the case when $\th$ is close to  $\frac12 \pi$  and the use of (\ref{K0}) does not lead to adequate result   we seek  a  finer estimation of the integral and to this end we   split the remaining  interval $[1/v,\infty)$ at  $v^{-1/4}$.  (For any number $\frac15<p< \frac13$, we may take  $ v^{-p}$ as the  
point of splitting  instead of $v^{-1/4}$.)

Put $$\a= \sqrt{1- \sin \th}$$ 
(so that  $|\frac12\pi -\th| \sim \sqrt{2}\,\a$).
\begin{Lem}\label{lem5.2.3} \quadd 
$J_{[v^{-1},  v^{-1/4}]} \leq  C/v\a^{3} \qquad \mbox{if}  \quad v \a^{3} \geq 1.$
\end{Lem}
\v2\n
\pf \,
  In place of Lemma \ref{lem3.4} we apply   Corollary \ref{prop_main} (in Section {\bf 4,2}), according to  which  we have
$$
  h_1(1+y,t, \phi) \leq \frac{Cy}{t^{2}} 
  \exp\Big\{- \frac{1}{2t} \eta^2(1 -  c\eta^4) \Big\}
  $$ 
with some universal constant $c$. In view of  the inequality 
 $\sqrt{\eta^2-    c \eta^6\,}> \eta (1- c\eta^4)$ ($0<\eta <\!<1$)  this application effects   replacing $K_1(v\sqrt{\phi^2+y^2})$ by 
 $ K_1(v\eta- cv \eta^5) $ 
  in the integral of (\ref{K0}), so that the exponent appearing in  its integrand  is at most 
   $$-v\sqrt{\eta^2- c\eta^6} +v\eta \sin \th \leq -v\a^2\eta + cv\eta^5.$$
     Hence    
$$J_{[v^{-1}, v^{-1/4}]}  \leq  C'  \int_{1/v}^{v^{-1/4}}  e^{- v\a^2 \eta + c v\eta^5 } \sqrt {v\eta }\,d\eta 
 $$
and the last integral is dominated by 
$$\frac{C'e^{c} }{v\a^3}\int_{\a^2}^{v^{3/4}\a^2}e^{-u}\sqrt u \, du \leq  \frac{C }{v\a^3},$$
as desired. \qed

\v2
\begin{Lem}\label{lem5.2.4}\quadd
$J_{[v^{-1/4}, \infty)} = O(ve^{-v^{1/12}}) \qquad \mbox{if}  \quad v \a^{3} \geq 1.$
\end{Lem}
\v2\n
\pf\,   
Lemma \ref{lem3.5} applied  with $t=s (<1)$ and  $r=\eta$ gives
\beqn\label{I4}
h_1^*(\xi^*_1 (l), s, \th) \leq \k_2 \frac{\eta}{s^2}\exp\Big\{-\frac{\eta^2}{2s}\Big\}. 
\eeqn
Substitution from  this bound in (\ref{K})  yields
\beqn\label{I41}
J_{[v^{-1/4}, \infty)} \leq C \int_{v^{-1/4}}^\infty e^{(1-\a)v\eta} d\eta \int_{0}^{1/v} \frac{\eta}{s^2}\exp\Big\{-\frac{v^2}{2}s -\frac{\eta^2}{2s} \Big\}ds. 
\eeqn
On applying (\ref{13}) again the inner integral on  the right-hand  side above  is  asymptotic to a constant multiple of 
$\sqrt{{v}/{\eta}}\,e^{-v\eta}$
as $v\eta\to \infty.$ Hence, for  $\a \geq v^{-1/3},$ 
\beq
J_{[v^{-1/4}, \infty)}  \leq C'\int_{v^{-1/4}}^\infty e^{-\a^2 v\eta}\sqrt{\frac{v}{\eta}} d\eta &=& \frac{C'}{\a }\int_{\a^2 v^{3/4}}^\infty e^{-u}\frac{d u}{\sqrt u} \\
&\leq& \frac{C''}{\a^2 v^{3/8}}e^{-\a^2 v^{3/4}} \leq C'''ve^{-v^{1/12}},
\eeq
where the last inequality follows from  $\a^2 v^{3/4} > (\a^2 v^{2/3})v^{1/12}$ and $\a^2 v^{3/8}> 1/v$.
Thus   the lemma  has been  proved. \qed
\v2

 Combining   Lemmas  \ref{lem5.2.2}, \ref{lem5.2.3} and \ref{lem5.2.4} we conclude (\ref{U0}) 
 as announced  at the beginning of Step 2.

\v2

 {\it Step 3.}  Here we compute  $U_{(a/v,t]}$ and finish the proof of Proposition \ref{UBD1}. We continue to suppose $a=1$.  Instead of (\ref{psi00}) we write 
 \beq
 e^{- x(1-\cos \th)/(t-s)}\, p^{(2)}_{t-s}(x -1) &=&  (1-s/t)^{-1}p^{(2)}_{t}(x-1) e^{- v(1-\cos \th)} \\
 &&\,\, \times \exp\Big\{\frac{- x^2s/t +2vs \cos \th - s/t}{2(t-s)}\Big\}
 \eeq
 and,  instead of (\ref{Ineq3}), we deduce the following expression of $\psi_1(l,t-s)e^{-\eta^2/2s}$:
\beq
&& \frac{x\cos\th -1}{t-s} \,[e^{- x(1-\cos \th)/(t-s)}p^{(2)}_{t-s}(x- 1) ]\,e^{x\eta (\sin \th)/(t-s)}e^{-\eta^2/2(t-s)}\times e^{-\eta^2/2s}\\
&&=\bigg(\frac{t}{t-s}\bigg)^{2}\frac{\Psi_1(x,t,\th)}{2\pi } 
 \exp\Big\{-\frac{1}{2(t-s)}\Big[\frac{x^2s}{t} +\frac{\eta^2t}{s} -2x\eta \sin \th - 2vs\cos \th +\frac{s}{t}\Big]\Big\}.
\eeq
Write  the  formula  in the square brackets in  the exponent  as   
$$\frac{t}{s}\Big(\frac{s}{t} x \sin \th -\eta\Big)^2 + s\Big[(x\cos \th -2)v\cos \th +\frac{1}{t} \,\Big]$$
and
apply Lemma \ref{lem3.5} to see  the bound 
$h_1^*(\xi^*_1(l), s, \th) \leq C_1  {\eta}{s^{-1}}p_s^{(2)}(\eta) $ for $s<1$. Then  we readily deduce that for $s<1$,
\beq 
\frac{\psi_1(l,t-s) h_1^*(\xi^*_1(l), s, \th)}{\Psi_1(x,t,\th)} &\leq& C\frac{t^2\eta}{(t-s)^2s^2} \exp\Big\{-\frac{t}{2(t-s)s}\Big(\frac{s}{t}x\sin \th -\eta \Big)^2\Big\} \\
&&\times  \,\exp \Big\{-\frac{s}{2(t-s)}\Big[(x \cos \th - 2)v\cos \th  \Big]\Big\}.
\eeq
We integrate the right-hand side over the half line $\eta\geq 0$. By applying  the inequality
$$\int_0^\infty p_T^{(1)}(\eta-m)\eta d\eta = \int_{-m}^\infty p_T^{(1)}(u) (u+m)du 
\leq \sqrt{\frac{T}{2\pi}}e^{-m^2/2T} + m$$
(valid for $m>0, T>0$), an easy computation yields 
$$\frac{U_{[1/v, 1/2]}}{\Psi_1(x,t,\th)}  \leq C'e^{-v/4}+
C' v\int_{1/v}^{1/2}  \exp \Big\{-\frac{s}{2(t-s)}\Big[(x \cos \th -2)v\cos \th  \Big]\Big\}\frac{ds}{\sqrt s}
$$
($v>2, t>1$),  of which the right-hand side is $O(e^{-\frac13 v^{1/3}})$ if $\cos \th \geq  v^{-1/3}$, hence $U_{[1/v, 1/2]}$ is negligible in this regime. We use Lemma \ref{lem3.5} again  to have  the bound 
$h_1^*(\xi^*_1(l), s, \th) \leq C_1 p_s^{(2)}(\eta) $ for $s\geq1/2$ and  we see $U_{[1/2, t]} = O(e^{-v})$ in a similar way. 

The  proof of Proposition  \ref{UBD1} is now complete.  \qed

\v2
The next lemma,  essentially a corollary of the proof of  Proposition \ref{UBD1},  provides a  crude upper bound for the case   $\cos \th \leq - v^{-1/3}$. Combined with Corollary \ref{cor1} it in particular verifies Corollary \ref{cor2}.
\begin{Lem} \label{wc}
$$h_a(x,t, \phi) \leq 
C \frac{-\Psi_a(x,t,\th)}{ |\th-\frac12 \pi|^3 v} \quad\quad \mbox{if}\quad \frac{\pi}2 +\frac1{(av)^{1/3}} <\th \leq \pi .
$$
\end{Lem}
\v2\n
\pf\, We have $h_a(x,t, \phi) \leq U$ (see (\ref{m9} and the  discussion succeeding it if necessary) and  observe  that  the identity (\ref{psi01}), hence the inequality (\ref{Ineq3}) are valid  for $\frac12 \pi <\th <\pi$ if the minus sign is put on the right-hand sides of them.  The proof of  (\ref{U}) may be then  adapted in a trivial way to the present case. \qed

\v2\n
{\bf 5.3.} {\sc Upper Bound II}. 
 
\begin{Prop}\label{UBD2} \,  Let $v=x/t >2/a$ and $t>a^2$. For some universal constant  $C$
$$h^*_a(\x,t,\th) \leq C (a v)^{2/3} e^{-av(1-\cos \th)} p_t^{(2)}(x -a) \quad \mbox{if} \quad \Big|\frac{\pi}2 -\th\Big| \leq  \frac1{(av)^{1/3}}.
$$
\end{Prop}
\v2\n
\pf\,      Let $a=1$. Put $\ga =\frac12 \pi -\th$ and suppose $ |\ga| \leq v^{-1/3}$. Let $\de$ be a small positive number chosen later and $\b =\ga+\de$ and  denote by  $L(\b)$  the line passing through the origin and $e^{i(\frac12\pi-\b)}$ so as to  make  the angle $\frac12 \pi-\b$ with the real axis (see Figure 2 in Section {\bf 5.4.} below). In this proof we consider the first hitting of $L(\b)$ by the two-dimensional Brownian motion starting at $\x= x$ (or $=x\be$). Let $y$ be the coordinate of $L(\b)$ such that $y=0$ for the point $e^{i(\frac12\pi-\b)}$ and $y=-1$ for the origin
and $\psi_\b(\x; y,t)$  the density of  the hitting distribution of $L(\b)$.  Let $\sigma(L(\b))$ denote the first hitting  time  of $L(\b)$ and $\eta(B_{\sigma(L(\b))})$   the $y$ coordinate of the hitting site  $B_{\sigma(L(\b))}\in L(\b)$. Then  we deduce 
\begin{eqnarray}\label{psi}    
\psi_\b(\x;y,t)& :=& \frac{P_\x[ \eta(B_{\sigma(L(\b))})\in dy, \sigma(L(\b)) \in dt] }{dydt}  \\
&=& \frac{x\cos \b}{t}p^{(1)}_{t}(x\cos \b)p^{(1)}_{t}(x\sin \b -y-1) \nonumber\\
&=&\frac{x\cos \b}{t}p_{t}^{(2)}(x)\exp\Big\{\frac{x(y+1)\sin \b -\frac12 (y+1)^2}{t}\Big\}.\nonumber
\end{eqnarray}
It  holds that
\beqn\label{clear}
h^*_1(\x,t,\th) \leq  2\int_0^t ds \int_{0}^\infty \psi_\b(\x;y,t-s)h^*_1(\xi^*(y),s,\th)dy,
\eeqn
where $\xi^*(y)$ denotes the point of $\R^2$ lying on $L(\b)$ of  coordinate $y$ (see Figure 2 of the next subsection). According to Lemmas \ref{lem3.4} and \ref{lem3.5}
\beqn\label{h_bd1}
h^*_1(\xi^*(y), s,\th) \leq \left\{
\begin{array}{ll} C ys^{-2} e^{-(y^2+\de^2)/2s} \quad &\mbox{if}  \,\,\, \,y<1, s<1,\\[2mm]
C(rs^{-1}\vee 1)p_s(r) \quad &\mbox{otherwise},
\end{array} \right.
\eeqn
where $r= |\xi^*(y) - e^{i\th}|$.
The rest of the proof is performed by showing  Lemmas \ref{lem_up_bd1}  and \ref{lem_up_bd2} given below.
\begin{Lem}\label{lem_up_bd1} \, For some universal constant $C$, 
 $$\int_0^{1/v} ds \int_{0}^{\infty} \psi_\b(\x; y,t-s)h^*_1(\xi^*(y),s,\th)dy \leq C v e^{v\cos \th} p_t^{(2)}(x) v^{-1/3}.$$
\end{Lem}
\v2\n
\pf\, We split the range of the outer integral at $y=1$ and  denote the corresponding the repeated integral  for  $[0,1]$ and $(1,\infty)$ by $I[0,1]$ and $I(1,\infty)$, respectively. 
As in the step 2 of the proof of Lemma \ref{UBD1} we see 
\begin{eqnarray}\label{**}
I_{[0,1]} &\leq& C vp_t^{(2)}(x)\int_{0}^1  e^{v(y+1)\sin \b} dy\int_0^{1/v} \frac{y}{s^2}e^{-\frac12 v^2 s - (y^2+\de^2)/2s} ds  \nonumber \\
&\leq&  C vp_t^{(2)}(x)\int_0^1  e^{v(y+1)\sin \b - v\sqrt{y^2+\de^2} } \frac{\sqrt v\,y}{(y^2+\de^2)^{3/4}}dy.
\end{eqnarray}
Put
$$f(y)= \frac{y^2}{\sqrt{y^2+\de^2} +\de} - 2\de y.$$
Suppose  $\de\geq \ga$. Then
\beq
\sqrt{y^2+\de^2} - (y+1)\sin \b &\geq& \sqrt{y^2+\de^2}-\de -2\de y -\sin \ga\\
&=& f(y) -\sin \ga
\eeq
and,  since $(x\sin \ga)/(t-s)= v\sin \ga + O(1)$ for $s\leq 1/v$ and 
 $\sin \ga =\cos \th$, 
 the last integral in (\ref{**}) is dominated by a constant multiple of 
\beq
&& e^{v\cos \th} \int_0^1 e^{-vf(y)}\frac{\sqrt v y}{(y^2+\de^2)^{3/4}}dy\\
&&= \frac{e^{v\cos \th}}{\sqrt{v\de}}\int_0^{ \sqrt{v/\de}} \exp\Big\{-\frac{u^2}{\sqrt{1+u^2/v\de}+1} +2\de^{3/2} v^{1/2}\, u\Big\}\frac{udu}{(1+ u^2/v\de)},
\eeq
where we have changed the variable of integration according to $y= (\de/v)^{1/2} \,u$. 
Now taking $\de =v^{-1/3}$ we can readily conclude that 
$$I_{[0,1]} \leq C v e^{v\cos \th} p_t^{(2)}(x) v^{-1/3}.$$

We can readily compute $I_{(1,\infty)}$  to be 
$v e^{v\cos \th} p_t^{(2)}(x) \times O(e^{-v/4}).$
Thus the proof of Lemma  \ref{lem_up_bd1} is complete. \qed
\begin{Lem}\label{lem_up_bd2} \, For some universal constant $C$, 
 $$\int_{1/v}^t ds \int_{0}^{\infty} \psi_\b(\x; y,t-s)h^*_1(\xi^*(y),s,\th)dy  \leq C v e^{v\cos \th} p_t^{(2)}(x) \times e^{-v/4}.$$
\end{Lem}
\v2\n
\pf\,  We  restrict the range of the outer integral  to $[1/v,1]$, the other part being easy to estimate, and   divide   the resulting  integral  by $p_t^{(2)}(x)$.
It suffices to examine the exponent of the exponential factor appearing  in $\psi_\b(\x;y,t-s)h^*_1(\xi^*(y),s,\th)/p_t^{(2)}(x) $ (in view of  (\ref{psi}) and (\ref{h_bd1}) ), which is 
\beq
&&-\,\frac{sx^2}{2t(t-s)}+ \frac{2x(y+1)\sin \b - (y+1)^2}{2(t-s)} - \frac{y^2+\de^2}{2s}\\
&&\leq -\, \frac{1}{2(t-s)}\Big(\frac{s}{t}x^2 + \frac{t}{s}y^2 -4x(y+1)\de\Big) - \frac{y}{2(t-s)}-\frac{\de^2}{2s}, 
\eeq
where for the inequality  we have applied $\sin \b \leq 2\de$. On the one hand for $y<1$, 
$$ [sx^2t^{-1}- 4x(y+1)\de]/2(t-s) \geq x(1-8\de)/2(t-s) \geq  v/3,$$
 provided $s\geq 1/v$ and  $\de<1/24$. On the other hand for $y\geq 1$
$$\Big(\frac{s}{t}x^2 + \frac{t}{s}y^2 -4x(y+1)\de\Big) = \bigg(\sqrt{\frac{s}{t}}x - \sqrt{\frac{t}{s}}y\bigg)^2
+2(1- 2(1+y^{-1})\de)xy,$$
which may be supposed larger than $vyt$. From these observations  it is easy to ascertain the bound of the lemma. \qed

 \begin{figure}[b]
 \begin{center}
\begin{picture}(150,105)

\put(116,24){\line (-5,1){100}}
\put(116,4){\line (-5,1){80}}
\put(36,20){\line (1,5){15}}
\put(51,-6){\line(1,5){20}}

\put(36,20){\circle{40}}
\put(36,20){\circle*{3}}
\put(40,39){\circle*{3}}
\put(59,35){\circle*{3}}

\put(58,18){$e^{-i\b}$}
\put(65,36){$e^{-i\b}(1+i)$}
\put(26,16){$0$}
\put(74,84){$l$}
\put(53,84){$y$}
\put(18,60){$L(\b)$}
\put(67,60){$L'(\b)$}
\put(40,-20){Figure 2}
\end{picture}
\end{center}
\vspace*{0cm}
\end{figure}

\v2\n
{\bf 5.4.} {\sc Lower Bound II and Completion of Proof of Theorem \ref{thm5.1}. }  

If $\cos \th > v^{-1/3}$,  the first  formula of Theorem \ref{thm5.1}  follows from Lemma \ref{LBD} and Proposition \ref{UBD1}.
Let  $|\cos \th| < v^{-1/3}$. The upper bound  in the second relation of Theorem  \ref{thm5.1}  follows from  Proposition \ref{UBD2}.  For derivation of  the lower bound we let $a=1$ and examine the proof of Lemma \ref{lem_up_bd1}. By the same computation as in it  with the help of the lower bound in Proposition \ref{cor_Imp} we see that 
$$I_{[\de^2, 1]} \geq  Cv e^{v\cos \th} p_t^{(2)}(x) v^{-1/3},$$ 
which however  is not enough since Brownian motion may have  hit $U(1)$ before $L(\b)$. 
The proof of the upper bound have rested  on the inequality (\ref{clear}), while  we  need a reverse inequality for the lower bound;  for the present purpose it suffices to prove
\beq
 h^*_1(\x, t,\th) \geq c \int_0^{1/v} ds \int_{\de^2}^1 \psi_\b(\x; y,t-s)h^*_1(\xi^*(y),s,\th)dy \eeq
for $|\frac12 \pi- \th|\leq v^{-1/3} $ and $\de=v^{-1/3}$ and for  some universal constant $c>0$, which, on comparing with (\ref{psi}),    follows if we have
\beqn\label{FP}
\psi_\b^*(\x; y,t) \geq c\psi_\b(\x; y,t) \quad \mbox{for} \quad \de^2< y<1
 \eeqn
(with the same  $c$ as above), where
$$\psi_\b^*(\x; y,t) =  \frac{P_\x[ \eta(B_{\sigma(L(\b))})\in dy, \sigma_{L(\b)} \in dt, \sigma_{1} >t] \,}{dydt} $$
($\eta(B_{\sigma(L(\b))})$ denotes the  $y$ coordinate of $B_{\sigma(L(\b))}$ as in the preceding proof). Let  $L'(\b)$ be the line tangent to the unit circle at  $e^{-i \b}$ and for the proof of (\ref{FP}) we consider the first hitting  by $B_t$ of  $L'(\b)$. Let $\z(l)$ denote the point on $L'(\b)$ of coordinate  $l$, where $l=0$ for $e^{-i\b}(1+i)$ and $l>0$ on the upper half of $L'(\b)$ (see Figure 2).  Then    for $\de^2<y<1$, we have
 \begin{eqnarray}\label{M} 
 \psi_\b(\x; y,t) &=& q_0^{(1)}(x\cos \b, t) p_t^{(1)}(y) \nonumber\\
 &=& \int_0^tds\int_{-\infty}^\infty q_{1}^{(1)}(x\cos \b, t-s)p_{t-s}^{(1)}(l) \psi_\b(\z(l); y, s)dl
 \end{eqnarray}
 and the corresponding  relation for $\psi_\b^*(\x; y,t)$ (with $\psi_\b^*$ in place of $\psi_\b$ in both places).
Noting $ \psi_\b(\z(l); y, s) = q_0^{(1)}(1,s) p_s^{(1)}(l-y)$ and integrating w.r.t.  $l$, we apply  Lemma \ref{9} (i)  given below (with $b=1$ so that $\rho t =1/v$ and $\sqrt{\rho t} =o(\de)$) (hence $(\rho t)^{3/2}  <\!< \de/v$)  to see that the outer integral may be restricted to  $|s-1/v| < \de/v$, so that 
$$\psi_\b(\x; y,t) 
\sim \int_{(1-\de)/v} ^{(1+\de)/v} ds\int_{-\infty}^\infty q_{1}^{(1)}(x\cos\b,t-s)p_{t-s}^{(1)}(l) \psi_\b(\z(l); y, s)dl,$$
of which  the inner integral  may be restricted to $l> 0$ with  at least  half the contribution of the integral preserved. Thus  the proof of  (\ref{FP}) is finished if we show that for some $c>0$,   $\psi_\b^*(\z(l); y, s) \geq c\psi_\b(\z(l); y, s)$ for $ y> s^{2/3}$ and $l\geq 0$, which we rewrite in terms of $\psi_{0}$ and $\psi^*_{0}$ as
\beqn\label{eq911}
\psi_0^*(1+i(1+l); y, s) \geq  c\,\psi_0(1+i(1+l); y, s) ,  \quad l\geq 0, \, y> s^{2/3}.
\eeqn
This is proved in Lemma \ref{91}   after showing the following  lemma.

\begin{Lem}\label{9}  Let  $0<b<x$ and  put  $\rho =b/x$.  For any $\e>0$ there exists a positive constant $M\geq 1$ that depends only on $\e$ such that {\rm (i)}  whenever $\rho t< 1/ M$,  $\rho <1-\e$ and $b\geq \e$,
\beqn\label{eq9}
\int_{|s- \rho t| < M(\rho t)^{3/2}} q_b^{(1)}(x, t-s) q_0^{(1)}(b, s)ds \geq (1-\e)q_0^{(1)}(x, t),\eeqn
and {\rm (ii)} whenever $ t< bx/ M^2$ and $\rho <1-\e$, (\ref{eq9}) holds if the range of integration is replaced by $|s - \rho t| < M(\rho t)^{3/2} b^{-1}$.
\end{Lem}

 The integral  in  (\ref{eq9}) extended to  the whole interval $[0,t]$  equals $q_0^{(1)}(x,t)$ and the lemma asserts that  substantial contribution to it  comes from a small interval about  $\rho t = bt/x$ (at least if  $x$ is kept  away from zero).
 \v2\n
\pf \, 
In this and the next proofs we apply the identity
  \beqn\label{eq90}
  p^{(1)}_{t-s}(z-\xi)p^{(1)}_s(y-z) =  p^{(1)}_{t}(y-\xi)p^{(1)}_{T}\Big( \frac{s}t (y-\xi) -y+z\Big),\quad T =\frac{s(t-s)}{t}
  \eeqn
  ($0<s<t, y,z, \xi\in \R$), This gives 
$$q_b^{(1)}(x, t-s) q_0^{(1)}(b, s) =\frac{(x-b)b}{(t-s)s} p^{(1)}_t(x)p^{(1)}_T\Big(\frac{s}tx -b\Big).$$ 
The range of integration of the integral in (\ref{eq9}) may be written as 
\beqn\label{eq903}
|s/t -\rho|\leq M\rho \sqrt {\rho t}, 
\eeqn
 which entails   $\frac{(x-b)b}{(t-s)s} =  \frac{(1-\rho)xb}{(1-s/t)ts} \sim\frac{xb}{ts}$ as $\rho t \to 0$,  and hence it suffices to show that  
\beqn\label{eq92}
\int_{|s- \rho t| < M(\rho t)^{3/2}} \frac{b}s \,p^{(1)}_T\Big(\frac{s}tx -b\Big)ds > 1-\frac12 \e
\eeqn
if   $1/\rho t$ and  $M$ are large enough.  Observing
$$\frac{b}s \,p^{(1)}_T\Big(\frac{s}tx -b\Big) =  \frac{b}{s\sqrt{2\pi (1-s/t)s\,}} \exp\Big\{ -\frac{b^2}{2(1-s/t)\rho t}\Big( \frac{s}{\rho t} +\frac{\rho t}{s} -2\Big)\Big\}$$
  and $u+u^{-1} -2 =  (1-u)^2 + O((1-u)^3)$ as $u\to 1$,
we apply the Laplace method  to see that  the integral in (\ref{eq92})  is asymptotic to  
 $$\int_{|u- 1| < M\sqrt{\rho t}}  \frac1{\sqrt{2\pi \la}}  e^{- (u-1)^2/2\la}du,$$
 where $\la = (1-\rho)\rho t/b^2$. If   the variable of integration is changed by $y = (u-1)/\sqrt \la$,
 then this integral  becomes $\int_{-r}^r p^{(1)}_1(y)dy$ with $r$ given by   
 $$r= M \sqrt{\rho t/\la} = Mb/\sqrt{1-\rho},$$
which extends to  the whole line as  $M\to\infty$ if $b\geq \e$.  Thus we obtain the assertion (i). 

As for the second assertion (ii)   we  multiply  the right-hand side of (\ref{eq903})  by   $b^{-1}$, and if  $b^{-1}\rho\sqrt{\rho t} =\sqrt{\rho t}/x = \sqrt{bt/x^3}\to 0$, then $\frac{(x-b)b}{(t-s)s}  \sim\frac{xb}{ts}$ as above. The rest of the proof  is the same.
\qed


\v2
Recall  (\ref{eq911}) and note that it  expresses the inequality
\beqn\label{eq91}
\frac{P_{1+i+il}[ \Im B_\tau -1 \in dy, \tau\in ds,  \tau <\sigma_1]}{ds dy} \geq  \frac{cP_{1+i+il}[ \Im B_\tau-1 \in dy, \tau\in ds]}{ds dy},
\eeqn
where  $B_t$ is a standard complex Brownian motion and  $\tau$ is the first hitting time of the imaginary axis by it.

\begin{Lem}\label{91} For a constant $c>0$,  (\ref{eq91}) holds true for $0< s\leq 1$, $l\geq 0$ and $y\geq s^{2/3}$.
\end{Lem}
\v2\n
\pf \,  The proof rests on the fact that if $Y_t $ denote the  linear Brownian motion $\Im B_t$,  then the conditional probability 
\beqn \label{923}
P[Y_{s'}>0, 0\leq s'\leq  s|Y_0=l, Y_s=y] = 1-e^{-2yl /s} \quad (l> 0, y>  0, s>0) 
\eeqn
 is bounded  away from zero if (and only if) so is $yl/s$. ((\ref{923}) is immediate from the expression of  transition density  for $Y_t$  killed at the origin.)

For  $\xi >0$ put
$$Q_\xi(y,t)=  q^{(1)}_0(\xi,t)p^{(1)}_t(y).$$
Then  for $0<b<1$,
\beq  
\psi_0(1+ i(1+l);y,s) &=&Q_{1}(y-l,s)\\
&=&  \int_0^s ds'\int_{-\infty}^\infty Q_{1-b}(y'-l, s-s') Q_{b}(y-y',s') dy'.
\eeq
Take $b= s^{1/3}$ in  the last integral.   Then by performing the integration w.r.t. $y'$ and noting  $(bs)^{3/2}b^{-1} = b^2s$ we apply   Lemma \ref{9} (ii)  (with $x=1$) to
infer  that  the $s'$-integration  above may be restricted to the interval 
$$ |s'-bs| \leq Mb^2s$$
with  some  $M\geq 1$.  Let $\phi= \tan^{-1} b$, $\eta= |b+i -e^{i\phi}| (= \sec\phi -1)$ and $\sigma(L_b)$ be the first hitting time of the line $L_b:=\{b+iy': y'\in \R\}$.   Since the slope of the tangent line of $\partial U(1)$ at $e^{i\phi}$ is  $b+o(b)$  and  
$$b\eta/s \sim 1/2$$
 (as $s\to 0$),  the  identity (\ref{923}) shows that if $s' \sim (1- b)s \sim s$,
$$\frac{P_{1+i(1+l)}[ \Im B_{\sigma(L_b)} \in dy', \sigma(L_b)\in ds', \sigma_1>s']}{dy'ds'} \geq c_1 Q_{1-b}(y',s'), \quad y'\geq 0$$
with  $c_1 =\frac12( 1- e^{-1})$, hence $\psi_0^*(1+i(1+l); y, s)$ is bounded below by a constant multiple of 
\[
  \int_{|s'-bs|<  Mb^2s}ds'\int_{y'\geq 0} Q_{1-b}( y'-l, s-s')\psi_0^*(b+i(1+y'); y, s') dy'.
\]
It therefore sufices to show that there exists $c_2>0$ such that if $s'\sim bs$ and $ y\geq s^{2/3}$, then
\beq
\psi_0^*(b+i(1+y'); y, s') \geq  c_2Q_{b}(y-y', s'),  \quad y'\geq 0, 
\eeq
which also follows from (\ref{923}) as is easily checked by noting   $s^{2/3} \eta/bs \sim \frac12$. Thus the lemma has been proved.
\qed
\v2


\section{Case $d\geq 3$ and Legendre Process}
This section consists of two subsections. 
The first one   concerns   the transition density of  a Legendre Process and  provides   the spectral expansion of it  as well as   its behavior  for small time which are  employed  in Section {\bf 3.3} and Section {\bf 4}, respectively. The second subsection  is devoted to the proof of Theorem \ref{thm1.2} for $d\geq 3$.

\v2\n
{{\bf 6.1.} \sc Legendre Process.} \,   

Let $d\geq 3$. The colatitude $\Th_t$ of $B_t/|B_t|$ is a Legendre process on $[0,\pi]$ regulated by the generator
$$\frac1{2\sin^{2\nu} \th}\frac{\partial}{\partial \th}\sin^{2\nu}\th\frac{\partial}{\partial \th}
=\frac12 \frac{\partial^2}{\partial \th^2} +\nu\cot \th\frac{\partial}{\partial \th}$$
with each boundary point  being entrance and non-exit (\cite{IM}). We compute the transition  law of $\Th_t$. Let $P_t^\nu(\th_0, \th)$ be the density of it
w. r. t. the normalized invariant measure: 
$$   \frac{P[\Th_t\in d\th|\Th_0=\th_0] }{d\th} =  P^\nu_t(\th_0, \th)\frac{\sin^{2\nu}\th}{\mu_d},$$
where $\mu_d=\int_0^\pi \sin^{d-2}\th d\th= \omega_{d-1}/\omega_{d-2}$. 

\v2
{\bf 6.1.1.} {\sc Eigenfunction Expansion.} \,
Eigenfunctions of the Legendre semigroup are  given by 
$$C_n^\nu(\cos \th) = \sum_{j=0}^n \frac{\Ga(\nu+j)\Ga(n+\nu- j)}{j!(n-j)![\Ga(\nu)]^2}\cos[(2j-n)\th],$$
where $C_n^\nu$ is a polynomial of order $n$ called the Gegenbauer (alias  ultraspherical)  polynomial   and in the special case $\nu=\frac12$ it agrees with the Legendre polynomial (see Appendix (A)).   
They together constitute a complete  orthogonal system of $L^2([0,\pi], \sin^{2\nu} \th d\th)$. (Cf.\cite{SW}, p.151 and 
\cite{W}, p.367; also \cite{T}, Section 4.5 for $\nu=1/2$.)  Given $\nu>0$, we denote their normalization by $h_n(\th)$:
$$h_n(\th)= \sqrt{\mu_d}\,\ga_n^{-1} C_n^\nu(\cos \th), $$
where  the factors $\ga_n>0$  are given by
$$\ga_n^2=\int_0^\pi  [C_n^\nu(\cos \th) ]^2\sin^{2\nu} \th d\th = \frac{\pi \Ga(n+2\nu)}{2^{2\nu-1} [\Ga(\nu)]^2(n+\nu)n!}$$
(cf. \cite{Sm}). (It is readily checked  that $\mu_d/\ga_0^2 =1$, so that $h_0\equiv1$.)
Then 
\beqn\label{spexp}
P^\nu_t(\th_0, \th) =  \sum_{n=0}^\infty e^{-\frac12 n(n+2\nu)t}h_n(\th_0)h_n(\th).
\eeqn
For translation of  the formula of Corollary \ref{cor-12} into that of Theorem \ref{thm1.21} one may use the formulae
 $C_n^\nu(1)= \Ga(n+2\nu)/\Ga(2\nu)n!$ and  $\Ga(2\nu) = 2^{2\nu-1}\Ga(\nu)\Ga(\nu+\frac12)/\sqrt \pi$ to see 
  $$h_n(0)h_n(\th) = \frac{\mu_d C_n^\nu(1)}{\ga_n^{2} }\, C_n^\nu(\cos \th) = \frac{\nu+n}{\nu}C_n^\nu(\cos \th) .$$
\v2
{\bf 6.1.2.} {\sc Evaluation of $P^\nu_t(0, \th)$ for $t$ small.} \, 
An  application of transformation of drift shows that uniformly for $ 0\leq \th <1$ and $t <1$
\beqn \label{asymp}
P^\nu_t(0, \th)  = \om_{d-1}p^{(d-1)}_t(\th) \Big[1+O(  \th^4 +t)\Big]. 
\eeqn
Indeed, if $X_t$ is a $(d-1)$-dimensional  Bessel process, $\ga(\th) =\nu( \cot \th - \th^{-1})$ and
$$Z_t = \exp\Big\{ \int_0^t\ga(X_s)dX_s -\int_0^t [\nu \ga(X_s)X_s^{-1} + \tst12 \ga^2(X_s)]ds\Big\}$$
then
$$P[\Th_t\in d\th,  {\cal E}_t^{\Th}\,|\,\Th_0=\th_0] = E^{BS(\nu-\frac12)}[ Z_t; X_t\in d\th,\, {\cal E}_t^{X} \,|\, X_0=0],$$
where ${\cal E}_t^{\Th}=\{\Th_s <1\,  \,\mbox{for}\,\, s<t\}$, ${\cal E}_t^{X}=\{X_s <1\,  \,\mbox{for}\,\, s<t\}$ and $E^{BS(\nu-\frac12)}$ signifies the expectation by the law of the Bessel process $X_t$. By simple computation using  Ito's formula we have  
$$Z_t= \exp\bigg\{\int_0^{X_t}\ga(u)du -\frac12 \int_0^t [\ga'(X_s)  +2\nu \ga(X_s)X_s^{-1}  + \ga^2(X_s)]ds\bigg\}$$ 
as well as $\ga(\th) = -\frac13 2\nu\th + O(\th^3), \ga'(\th)= -\frac13 2\nu +O(\th^2)$. 
Noting that  $p_t^{(d-1)}(\th)$ is the density of $P^{BS(\nu-\frac12)} [X_t \in d\th \,|\,X_0=0]$ w.r.t.   
 $\om_{d-2}\th^{d-2}d\th$ and 
$$[\om_{d-2}\th^{d-2}]/ [\mu_d^{-1}\sin^{2\nu} \th] = \om_{d-1}(1 + 3^{-1} \nu \th^2) + O(\th^4),$$  substitution yields (\ref{asymp}). 

\v2\n
{\bf 6.2.} {\sc Proof of Theorem \ref{thm1.2} ($d\geq 3$).}  

Recall  the definitions of     $g(\th;x,t)$ given  in Section {\bf 3.3.1}  and of
$h_a( x,t, \phi)$   in (\ref{h_z00}). Noting $|d\xi|=a^{d-1}\om_{d-1} m_a(d\xi)$, we then see that for  $\x=x\be$ and $\xi \in \partial U(a)$ of colatitude $\th$ 
$$\frac{P_{\x}[B_{\sigma_a}\in d\xi, \sigma_a\in dt]}{|d\xi|dt} = \frac{g(\th;x,t)}{a^{d-1}\om_{d-1}} q^{(d)}(x,t)=   \frac{h_a(x,t, \th)}{a^{d-1}\om_{d-1}}
$$
and  that owing to Theorem A  the two relations of  Theorem \ref{thm1.2}  are  equivalent to the corresponding ones   in Theorem \ref{thm5.1}  if adapted to the higher dimensions:   in the right-hand side of the first formula of   Theorem \ref{thm5.1}    the heading factor $2\pi a$
is replaced by $a^{d-1}\om_{d-1}$ and $p_t^{(2)}(x-a)$ by $p_t^{(d)}(x-a)$, and similarly for  the second one.
For the proof of them
 we may repeat the same procedure for two-dimensional case  with suitable modification, but here  we adopt another way  of reducing the problem to that for the two-dimensional case: roughly speaking we 
have   $(d-2)$-dimensional  variable, denoted  by  $\z$,  against  which the additional factor 
\beqn\label{factor}
p^{(d-2)}_{t-s}(z)p^{(d-2)}_s(z), \quad z=|\z|
\eeqn
that must be incorporated in the computation is integrated to produce the factor $p_t^{(d-2)}(0) $ (because of  the semi-group property of $p_t$),  which together with $p_t^{(2)}(x-a)$    constitutes 
the factor $p_t^{(d)}(x-a)$ in the final formula.

More details are given below. Recollecting the proof of Proposition \ref{UBD1}, we regard  the two-dimensional space where the problem is discussed in it  as a subspace of $\R^d$ in this proof and the line $L(\th)$ (introduced in the proof of Lemma \ref{LBD}) as the intersection of this subspace with a $(d-1)$-dimensional hyper-plane, named $\De(\th)$,  that is tangent at  $\xi$ with $\xi\cdot \be =\cos \th$  to the sphere $\partial U(a)$. (Here we write  $\De(\th)$  for the hyper-plane which is determined not by  $\th$ but by $\xi$ since the variable   $\th$ is essential in the present  issue.)  Let $M(\th, l)$ be the $(d-2)$-dimensional subspace contained in $\De(\th)$ passing through  $\xi^*(l)\in L(\th)$ ($l$ is a coordinate of $L(\th)$ as before) and perpendicular to the line $L(\th)$. Put
$$H_a(\y, t,\xi) = \frac{P_\y[B(\sigma_a)\in d\xi, \sigma_a\in dt]\,}{m_a(d\xi)dt}\quad\quad  (\y \notin U(a),  \,\xi\in \partial U(a)) $$
and
$$ \psi(l,t) = \frac{P_\x[\, {\rm pr}_{L(\th)} B(\sigma_{\De(\th)})\in d l, \sigma_{\De(\th)}\in dt]}{dl dt} ,$$
where  $ {\rm pr}_{L(\th)}$ denotes the orthogonal projection  on $L(\th)$, and  define $U^{(d)}$  as in (\ref{K}) but with $H_a(\y,  t,\xi)$ in place of $h^*_a(\y, t,\th)$.
 Then for each $l$ the claim (\ref{U}) is replaced by 
\beq
U^{(d)}  &=& \int_0^t ds \int_{\R} \psi(l, t-s)dl \int_{M(\th,l)} p^{(d-2)}_{t-s}(z)H_a(\xi^*(l) +\z,s,\xi) |d\z|   \\
&\leq & \frac{C\Psi_a(x,t,\th)}{av\cos^3 \th}.
\eeq

For the region in which  $z < \eta $, we may simply multiply the integrand in (\ref{K}) by  (\ref{factor}) without anything that requires  particular attention. If  $z > \eta $, we also multiply  the integrand by  $p^{(d-2)}_{t-s}(\z)$, replace $h^*_a(\xi^*(l), s,\th)$ by $H_a(\xi^*(l)+\z,s, \xi)$  and  use the bound 
$$ H_a(\xi^*(l)+\z,s,\xi) \leq  \frac{C z^2}{s}p_s^{(2)}(\eta) p_s^{(d-2)}(z)\Big(1\vee\frac{z^2}{\sqrt s}\Big)  e^{C_1z^6/2s}$$
in  Step 2 (Lemmas \ref{lem5.2.2} and  \ref{lem5.2.4}) (i.e., the step corresponding to that in the proof of Proposition \ref{UBD1}), and 
$$ H_a(\xi^*(l)+\z,s,\xi) \leq  \frac{C z}{s}p_s^{(2)}(\eta) p_s^{(d-2)}(z)$$
in the last part of Step 2 (Lemma \ref{lem5.2.4}) and in Step  3.  In  Step 2  there appears the integral
$$\int_0^b \bigg(\frac{z^2}{s} + \frac{z^4}{\sqrt s}\bigg) \exp\Big\{-\frac{z^2 -6C_1z^6}{2T}\Big\} \frac{z^{d-3}dz}{T^{(d-2)/2}} \quad \mbox{where}\quad T= \frac{s(t-s)}{t},$$
which is made less than unity   for $s<1/v$ by taking $b$ small enough, especially with  $b=v^{-1/4}$.
 In  Step 3 (and  the last part of Step 2) we have only to notice that
 $$\int_b^\infty \frac{z}{\sqrt s} p^{(d-2)}_{t-s}(z)p^{(d-2)}_s(z)z^{d-3}dz \leq C e^{-vb/4} $$
for $s<1/v$. 
With these considerations  taken into account  the proof of Proposition \ref{UBD1} goes through  virtually intact. The
further details are omitted. 

In a similar way  Proposition  \ref{UBD2} and the lower bound obtained in  {\bf 5.3} are extended to the dimensions $d\geq 3$.

 \v2\v2\n
\section{ Brownian Motion with A Constant Drift} 
In this section we present the results for the  Brownian motion with a constant drift  that are readily derived from those given above for the bridge.
The Brownian motion  $B_t$  started at $\x$ and   conditioned to hit $U(a)$ at $t$  with $v:=x/t$ kept away from zero may be comparable or similar to  the process $B_t - tv{\bf e}$ in significant  respects and  some of our results for the former one is more naturally comprehensible   in its translation in terms of the latter (see (\ref{transl}) at the end of this section).

\v2\n
{\bf 7.1. \, Formulae in general setting}
\v2\n
Given $v>0$, we put  
$$\bv = v\be$$
 (but $\x\notin U(a)$ is arbitrary) and label the objects defined by means of   $B^{(\bv)}_t:=B_t- t {\bf v}$  in place of $B_t$ with the superscript  $\,^{(\bv)}$ like $\sigma_a^{(\bv)}, \Th_t^{(\bv)},$ etc.  The translation is made by using the formula for  drift transform.  We put   $\ga(\cdot)= -\bv$ (constant function) and
$Z(s) = e^{\int_0^s \ga(B_u)\cdot d B_u -\frac12 \int_0^s |\ga|^2(B_u)du}$, so that $P_\x[ (B_t^{(\bv)})_{t\leq s}\in \Ga] = E_\x[ Z(s); (B_t)_{t\leq s}  \in \Ga]$ for $\Ga $ a measurable set of $C([0,s],\R^d)$.
 It follows that    $Z(\sigma_a) = \exp\{-\bv\cdot B(\sigma_a) +\bv\cdot B_0-\frac12 v^2\sigma_a \}$. Hence   
 \beq
 &&P_{\x}[B^{(\bv)}(\sigma^{(\bv)}_a) \in d\xi, \sigma^{(\bv)}_a \in dt] \\
 &&\quad = e^{\bv\cdot\x  -\frac12 v^2t}e^{-\bv\cdot \xi}P_{\x}[B(\sigma_a)\in d\xi, \sigma_a \in dt],
 \eeq
 and
putting
 $$  f_{a,t}^{(\bv)}(\x,\xi)   =\frac{e^{-\bv\cdot \xi}P_{\x}[ B(\sigma_a) \in d\xi \,|\,  \sigma_a =t]}{m_a(d\xi)},
$$
we obtain
$$
\frac{P_{\x}[B^{(\bv)}(\sigma^{(\bv)}_a) \in d\xi, \sigma^{(\bv)}_a \in dt] }{m_a(d\xi)dt}=  e^{\bv\cdot\x  -\frac12 v^2t}q_a^{(d)}(x,t) f_{a,t}^{(\bv)}(\x,\xi),
$$
\beqn\label{7.1.1}
\frac{P_{\x}[ \sigma^{(\bv)}_a \in dt] }{dt}=  e^{\bv\cdot\x  -\frac12 v^2t}q_a^{(d)}(x,t) \int_{\partial U(a)} f_{a,t}^{(\bv)}(\x,\xi)m_a(d\xi),
\eeqn
and
\beqn\label{7.1.2}
 \frac{P_{\x}[B^{(\bv)}(\sigma^{(\bv)}_a) \in d\xi \,|\, \sigma^{(\bv)}_a  = t] }{m_a(d\xi)} =  \frac{f_{a,t}^{(\bv)}(\x,\xi)}{\int_{\partial U(a)} f_{a,t}^{(\bv)}(\x,\xi)m_a(d\xi) }.
\eeqn


Suppose  $x/t\to 0$ and $t\to\infty$.  By Theorem \ref{thm1.1}, 
$$f_{a,t}^{(\bv)}(\x,\xi) =e^{-\bv\cdot \xi} \Big(1+ O\Big(\frac{x}{t}\ell(x,t)\Big)\Big),$$
so that 
$$\frac{P_{\x}[\sigma^{(\bv)}_a \in dt] }{dt} = \Big[\int_{|\xi|=a}  e^{-\bv\cdot \xi}m_a(d\xi)\Big] e^{\bv\cdot\x  -\frac12 v^2t}q_a^{(d)}(x,t)  \Big(1+ O\Big(\frac{x}{t}\ell(x,t)\Big)\Big),$$
where  $\ell(x,t)$   is the same function as given in Theorem \ref{thm1.1} if $d=2$ and $\ell(x,t) \equiv 1$  if $d\geq 3$.
Noting  $e^{\bv\cdot \x -\frac12 v^2t}p_t^{(d)}(x) = p_t^{(d)}(|\x - t\bv|)$  we deduce from  Theorem A  that
\beqn\label{7.1.3}
e^{\bv\cdot\x  -\frac12 v^2t}q_a^{(d)}(x,t) = a^{2\nu}\La_\nu\Big(\frac{ax}{t}\Big)p_t^{(d)}(|\x - t\bv|) \bigg[1-\Big(\frac{a}x\Big)^{2\nu}\bigg](1+o(1))
\eeqn
for $d\geq 3$ and  an analogous relation for $d=2$ (where  the formula must be modified  in the case $x\leq \sqrt t$
according to (\ref{R21})). We have   the identities  $C_0^\nu \equiv 1$ and 
$$
\int_0^\pi e^{-z\cos \th}  C_n^{\nu}(\cos \th) \sin^{2\nu}\th \,d\th 
=(-1)^n\frac{2^\nu \sqrt \pi \Ga(\nu+\frac12)\Ga(n+2\nu)}{\Ga(2\nu)n!} \cdot \frac{I_{n+\nu}(z)}{z^\nu},$$
where $I_\nu(z)$ is the modified Bessel function of the first kind of order  $\nu$ and, on putting $n=0$ in the latter, 
$$\int_{|\xi|=a} e^{-\bv\cdot \xi}m_a(d\xi) = \frac{2^\nu \sqrt \pi \,\Ga(\nu+\frac12) }{\mu_{d}}\cdot\frac{I_{\nu}(v)}{v^\nu}.$$

Let  $g(\phi;y)$, $y>0$,  denote the  function represented by the series in (\ref{lim}), namely
\beqn\label{g7.1}
g(\phi; y) = \sum_{n=0}^\infty  \frac{K_\nu(y)}{K_{\nu+n}(y)}H_n(\phi).
\eeqn
Then, owing to Theorem \ref{thm1.21}, as  $x/t \to \tilde v>0$ and $t\to\infty$,
$$f_{a,t}^{(\bv)}(\x,\xi) =  e^{-\bv\cdot \xi} g(\phi ;a\tilde v)(1+o(1)) \quad \mbox{for}\quad \xi\in \partial U(a),
$$
where $\xi\cdot\x/ax = \cos \phi$. It is worth noting that  if $ \tilde v/v$ is small, then the  function $e^{-\bv\cdot \xi} g(\th ;a\tilde v)$  is maximized about $\xi_a := - a\be\in  \partial U(a)$ (not $a\be$) irrespective of $\x$.
\v2\n
{\bf 7.2.  \, Case $\x- t\bv =o(t)$}

\v2\n
 In this subsection we  let $\x = x\be $,  while  $\bv$ is arbitrary but subject to the condition
$$\quad  \frac{ \x}{t}-\bv \,\longrightarrow \,  0,$$
 so that
$$\bv\cdot\xi = t^{-1}\x\cdot \xi +o(1) =  av\cos \th +o(1)$$
uniformly  for  $\xi\in \partial U(a)$ with $\xi\cdot \x/ax =\cos \th$.
  Define $g^{(\bv)}_a(\x,t,\th)$ by 
\[
g^{(\bv)}_a(\x,t,\th) =   \frac{P_{\x}[B^{(\bv)}(\sigma^{(\bv)}_a) \in d\xi\,|\, \sigma^{(\bv)}_a = t] }{m_a(d\xi)}.
 \]
Then by  (\ref{7.1.2})
\begin{equation}\label{7.2.2}
g^{(\bv)}_a(\x,t,\th) =
 \frac{e^{-\bv\cdot\xi}g_a(x,t,\th)}
{\mu_d^{-1} \int_0^\pi e^{-\bv\cdot\xi}g_a(x,t,\phi)\sin^{d-2}\phi \,d\phi},
\end{equation}
where $g_a(x,t,\th$ is defined in (\ref{g7.1}).
Let $\Xi_{av} $ denote the (normalizing) constant
$$\Xi_{av}=\int_{0}^\pi e^{- av\cos \th}g(\th;av)\frac{\sin^{d-2} \th  \, d\th}{\mu_d}.$$
(Remember that $g(\th;av)$   is the  density w.r.t. $\mu_d^{-1}\sin^{d-2} \th d\th$ of the limit distribution of $\Th(\sigma_a)$ conditioned on  $\sigma_a =t$,  $B_0=x\be$.)  Then as   $t\to \infty$ under  $|\x/t -\bv| \to 0$,
we have $\Xi_{av} \sim  \int_{|\xi|=a} f_{a,t}^{(\bv)}(\x,\xi)m_a(d\xi)$ and hence
\vskip3mm
 (i)   \quad ${\displaystyle   \frac{P_{x\be}[ \sigma^{(\bv)}_a \in dt] }{dt} \, \sim 
 \,\Xi_{av}\,a^{2\nu}\La_\nu(av) p_t^{(d)}(|\x - t\bv|) (1+o(1)); } $
 \vskip3mm
(ii) \quad  ${\displaystyle  g^{(\bv)}_a(x\be, t,\th) \,\sim\,\frac1{\Xi_{av}} e^{-av\cos \th}g(\th;av)}, $
\v2\n
where the last asymptotic relation  is uniform for $0\leq \th\leq \pi$ and $v \leq  M $ for any $M>0$; for (i)  use   the identities (\ref{7.1.1})  and  (\ref{7.1.3}). 
\vskip3mm\n

Similarly, substituting  the formula of Corollary \ref{thm3.02} in (\ref{7.2.2})   (cf.  (\ref{sin}))  we obtain 
an asymptotic form of $g_a^{(\bv)}(\x,t,\th)$ as $v\to\infty$.  On observing that this leads to 
   $$\Xi_{av} \sim \mu_{d}^{-1}v\int_0^{\pi/2}\sin^{d-2}\th \cos  \th d\th = \om_{d-2}v/(d-1)\om_{d-1}$$
    ($v\to \infty)$,  a simple computation  yields  the following asymptotic relations:  as $ v\to\infty$ and $|x -tv|/t \to 0$,
$$ \frac{P_{x\be}[ \sigma^{(\bv)}_a \in dt] }{dt}  \, \sim \, \frac{\om_{d-2}}{d-1}a^{2\nu+1}v\, p_t^{(d)}(|\x - t\bv|)$$
and, if $(av)^{-1/3} \leq \cos\th \leq 1$,
\beq
g^{(\bv)}_a(x\be, t,\th) 
  = (d-1)\mu_d \bigg[\cos \th+ O\Big(\frac{1}{av\cos^2 \th}\Big)\bigg](1+o(1)), 
\eeq
 where $o(1)$ is independent of  $\th$;
 also  Corollary \ref{thm3.01} may  translate into 
\beqn\label{E}
P_{x \bf e} [\Th_t^{(\bv)}\in d\th\,|\, \sigma_a^{(\bv)}=t]     \,\Longrightarrow \,(d-1)\1(0\leq\th<\tst12 \pi) \sin^{d-2}\th\, \cos \th\, {d\th}.
\eeqn
The last  convergence result  may be  intuitively comprehended  by noticing that  the right-hand side is the law  of  the colatitude of a random variable taking values in  the \lq northern  hemisphere' of  $\partial U(a)$ whose projection on the \lq equatorial  plane'    is uniformly distributed on the \lq\lq hyper  disc'', $\mathbb{D}$ say, on the plane; in short it may be thought as the distribution on the sphere induced  by the uniform ray coming from the direction $\be$. Let ${\rm pr}_{\be}$ denote this projection on the equatorial  plane. Then the  result given in  (\ref{E}) may be restated as follows:   $P_{x\bf e} [\pr_\be B^{(\bv)}_t\in dw\,|\, \sigma_a^{(\bv)}=t]$, $d w \subset \mathbb D$ converges weakly to the uniform measure on  $\mathbb{D}$. We rephrase Theorem \ref{thm1.2} in a similar fashion.  Let  $\xi \in \partial U(a)$,  $\xi\cdot \be =a\cos \th$ and $w={\rm pr}_{\be} \xi$ and note that
$$a-|w|  \sim 2^{-1}a\cos^2 \th \,\,\,  (\th \to {\textstyle \frac12} \pi), \quad  \cos \th = \sqrt{1-|w|^2/a^2} \quad\mbox{ and} \quad  |d\xi| = |dw|/\cos \th$$
and  that $m_a(d\xi) = a^{-d+1}|d\xi|/\om_{d-1}$  and $\om_{d-2}=(d-1)c_{d-1}^*$, where  $c^*_{n}$ denotes the volume of the  unit ball in  $\R^n$.   
Then, from Theorem \ref{thm1.2} we deduce that  uniformly  for   $w\in \mathbb{D}$,
\begin{eqnarray} \label{transl}
&& \frac{P_{x\bf e} [\pr_\be B^{(\bv)}_t\in dw\,|\, \sigma_a^{(\bv)}=t] }{[a^{d-1}c^*_{d-1}]^{-1}|dw|}  \\
 &&\quad = \bigg[1 +O\bigg(\frac1{(1-|w|/a)^{3/2} av}\bigg) \bigg](1+o(1))   \qquad \mbox{for} \quad   |w|/a <1- (av)^{-2/3},\nonumber\\
&& \quad \asymp    (av)^{-1/3}/\sqrt{1- |w|/a} \quadd \quadd \mbox{for} \quad  1- (av)^{-2/3} \leq |w|/a \leq 1, \nonumber
\end{eqnarray}
as $v\to \infty$ and $|x\be/t-\bv|\to 0$, showing  convergence of the density on the  one hand 
and indicating the effect of Brownian noise that  manifests  itself as the singularity of the density  along the boundary of $\mathbb{D}$.

The strict  equalities ${\bf x}/x ={\bf v}/v = {\bf e}$ we have assumed above can be relaxed. Essential assumption is
 ${\bf x} - t{\bf v} =o(t)$, entailing that  ${\bf v}\cdot \xi = t^{-1}{\bf x}\cdot \xi +o(1) =  av\cos \theta +o(1)$
uniformly  for  $\xi\in \partial U(a)$ with $\xi\cdot {\bf x}/ax =\cos \theta$. The identity (7.1) does not hold any more, but two  sides of it are asymptotically equivalent and the other relations including  (7.2) remain valid.


\section{Appendix}
(A) \, The Gegenbauer polynomials   $C^\nu_n(x)$, $n=0, 1,2,\ldots$,  may be   defined as the coefficients  of  $z^n$ in the Taylor series   $(z^2 -  2xz +1)^{-\nu}= \sum C^\nu_n(x)z^n$ ($|z|<1, |x|\leq 1, \nu>0$) and form an orthogonal basis of the space $L^2([-1,1],(1-x)^{\nu})$ (cf. page 151 of \cite{SW}).  
 The function $u(x)=C^\nu_n(x)$  satisfies
$$(x^2-1)u'' + (2\nu+1) xu' - n(n+2\nu)u =0
$$
and it follows that if  $Y(\th)=u(\cos \th)$,
$$\frac12 Y'' + \nu \cot \th\, Y' + \frac{n(n+\nu)}{2} Y =0.$$ 
\v2\n
 (B) The density $ P_\x[B(\sigma_a)\in d\xi, \sigma_a\in dt]/ m_a(d\xi)dt$ admits  an  explicit  eigenfunction expansion. In the case 
 $d=2$ it is given below. Let $p^0_{(a)}(t,\x,\y)$ denote the transition probability  of a two-dimensional Brownian motion that is killed when it hits $U(a)$. Then according to Eq(8) on p. 378 in \cite{CJ}
\beqn\label{seriesexp}
p^0_{(a)}(t, x\be, \y) = \frac{1}{2\pi}  \sum_{n=-\infty}^\infty \cos n\th 
 \int_0^\infty   e^{-\la^2 t/2} \frac{U_n(\la, x)U_n(\la, y)}{J^2_n(a\la)+ Y^2_n(a\la)}\la d\la,
 \eeqn
where  $J_n$ and  $Y_n$ are  the usual Bessel functions   of the first and second kind, respectively,
$$U_n(\la, y) =  Y_n(\la a) J_n(\la y) -  J_n(\la a)Y_n(\la y)$$
and $\y =(y,\th)$, the polar coordinate of $\y$ (with $y=|\y|$, $\cos \th =\y\cdot \be/y$).   
From the identity $(Y_\nu J'_\nu - J_\nu Y'_\nu)(z) = -2/\pi z$  it follows that $(\partial/\partial y) U_n(\la, y)|_{y=a} = -2/\pi a$ and 
\beq
\frac{P_{x\be}[{\rm Arg} \,B_t\in d\th, \sigma_a \in dt]}{a d\th dt} &=& \frac12\frac{\partial}{\partial y}  p^0_{(a)}(t, x\be, \y)|_{y=a} \\
& =&  \sum_{n=-\infty}^\infty  I_n(x, t)\cos n\th
\eeq
where
$$ I_n(x, t)= \frac{1}{2a\pi^2}  
 \int_0^\infty   e^{-\la^2 t/2} \frac{- U_n(\la, x) \la d\la }{J^2_n(a\la)+ Y^2_n(a\la)}.$$
Since integration  by $a d\th$ reduces   the density given above  to $q^{(2)}_a(x,t)$, we have
$$q^{(2)}_a(x,t) = 2\pi  aI_0(x, t) = \frac1{\pi} \int_0^\infty   e^{-\la^2 t/2} \frac{- U_0(\la,  x) \la d\la }{J^2_0(a\la)+ Y^2_0(a\la)}$$
and
$$\frac{P_{x\be}[{\rm Arg}\, B_t\in d\th \,|\,\sigma_a = t]}{2\pi d\th}= \frac1{2\pi}+\frac1{2\pi I_0(x,t)} \sum_{n=1}^\infty  I_n(x, t)\cos n\th.
$$
On comparing with (\ref{-0})  $2I_n(x,t)$ must agree with $a^{-1}q^{(2)}_a(x,t)\a_n(x,t)$, so that  
$$q_a^{(2n+2)}(x,t) = \bigg(\frac{a}{x}\bigg)^{n} 2a\pi I_{n}(x,t)
=\frac{1}{\pi}  \bigg(\frac{a}{x}\bigg)^{n} 
 \int_0^\infty   \frac{- U_n(\sqrt{2\a}, x)  e^{-\a t} d\a }{J^2_n(a\sqrt{2\a})+ Y^2_n(a \sqrt{2\a})}.$$
This last  formula, though not used in this paper, is valid for non-integral $n$ and  useful: e.g.,  its use  provides another approach in which one may dispense with  the arguments  using the Cauchy integral theorem for  the proofs   in \cite{Ubh} and \cite{Ubes}. 

The integral transform involved in the Fourier series  (\ref{seriesexp}) is derived by using the Weber formula (\cite{T1}, p. 86)  and  the higher-dimensional analogue  is given by   the Legendre series (as in (\ref{spexp}) with  
an  integral transform similar to the one in (\ref{seriesexp})). 

\v2\n
(C)\, We prove that  for each $\de>0$, as $y\downarrow 0$ and  $\phi \to 0$
\beqn\label{Rem510}
 \frac1{a^{d-1}\om_{d-1}}\int_0^\de h_a(a+y, s, \phi)ds = \frac{2y}{\om_{d-1} [y^2+ (a\phi)^2]^{d/2}}(1+o(1))
 \eeqn
(a result used in Remark 5). This is an expression of  the obvious fact that as  $y \downarrow 0$   the hitting distribution of $\partial U(a)$  for the Brownian motion started at 
$(a+y)\be$  converges to   that of the plane tangent to it at $a\be$: the ratio on the right-hand side is a substitute of the density of the latter distribution, where   $(a\phi)^2$ in the denominator must be replaced by  $|\z- a\be|^2$ with   $\z$ being  any point of the plane such that  $\z\cdot \be/ |\z| = \cos \th$.   

For verification let $P(\z,\xi;a)$ be the Poisson kernel of the exterior of the ball $U(a)$,  with respect to the 
 uniform probability $m_a(d\xi)$ so that $\int_{\partial U(a)}P(\z,\xi; a)m_a(d\xi) = (a/z)^{2\nu}$  ($z=|\z|>a$). It is given by
 $$P(\z,\xi; a) = \frac{a^{2\nu}(z^2-a^2)}{|\z-\xi|^d}, \quad  z>a, \, \xi\in \partial U(a).$$
Let  $\xi$ be such that  $\z\cdot\xi/za=\cos \phi$. Then by  an elementary  computation we find that
$$P(\z,\xi; a) = \frac{2a^{2\nu+1}y}{[y^2+ (a\phi)^2]^{d/2}}(1+o(1))\quad \mbox{as}\quad y:= z-a \downarrow 0, \,\,
\phi \to 0$$
and this shows (\ref{Rem510}), for    $P(\z,\xi; a)$ equals 
the whole integral $\int_0^\infty h_a(\z, s, \phi)ds$ and  this integral restricted  to  $[\de,\infty)$ is dominated by a constant multiple of $y$ owing to Theorem A and Theorem \ref{thm1.1}  (cf.  the first inequality  of Lemma \ref{lem3.5}). 

\v2

\end{document}